\documentclass[10pt]{article}
\usepackage{verbatim,latexsym,amssymb,amsthm,amsmath,amsfonts}

\def\R{{\hbox{\bf R}}}
\def\C{{\hbox{\bf C}}}

\def\B{{{\mathcal B}}}
\def\D{{{\mathcal D}}}
\def\E{{\hbox{\bf E}}}
\def\Energy{{{\mathcal E}}}

\def\P{{\hbox{\bf P}}}

\def\dist{{\hbox{\rm dist}}}
\def\mod{{\ \hbox{\rm mod}\ }}
\def\c{{\hbox{\bf c}}}
 at 10 true pt
\def\Z{{\hbox{\bf Z}}}
\def\eps{\varepsilon}

\def\emph#1{{\it #1}}
\def\textbf#1{{\bf #1}}
\newcommand{\vdW}{{\mbox{\scriptsize vdW}}}
\newcommand{\FAN}{{\mbox{\scriptsize FAN}}}
\newcommand{\SZ}{{\mbox{\scriptsize SZ}}}

\theoremstyle{plain}
  \newtheorem{theorem}[subsection]{Theorem}

  \newtheorem{proposition}[subsection]{Proposition}
  
  \newtheorem{lemma}[subsection]{Lemma}

\theoremstyle{remark}
  \newtheorem{remark}[subsection]{Remark}
  
  \newtheorem{example}[subsection]{Example}

\theoremstyle{definition}
  \newtheorem{definition}[subsection]{Definition}

\include{psfig}

\begin{document}

\title{A quantitative ergodic theory proof of Szemer\'edi's theorem}

\author{Terence Tao \\
\small Department of Mathematics, UCLA, Los Angeles CA 90095-1555 \\[-0.8ex]
\small tao@math.ucla.edu \\[-0.8ex]
\small http://www.math.ucla.edu/$\sim$tao \\[-0.8ex] }

\date{\small MR subject classification: 11B25, 37A45}

\begin{abstract}  A famous theorem of Szemer\'edi asserts that given any density $0 < \delta \leq 1$ and any 
integer $k \geq 3$, 
any set of integers with density $\delta$ will contain infinitely many proper arithmetic progressions of length $k$.
For general $k$ there are essentially four known proofs of this fact; Szemer\'edi's original combinatorial proof using the Szemer\'edi regularity lemma and van der Waerden's theorem, Furstenberg's proof using ergodic theory, Gowers' proof using Fourier analysis and the inverse theory of additive combinatorics, and the more recent proofs of Gowers and R\"odl-Skokan using a hypergraph regularity lemma.  Of these four, the ergodic theory proof is arguably the shortest, but also the least elementary, requiring in particular the use of transfinite induction (and thus the axiom of choice), decomposing a general ergodic system as the weakly mixing extension of a transfinite tower of compact extensions.  Here we present a quantitative, self-contained version of this ergodic theory proof, and which is ``elementary'' in the sense that it does not require the axiom of choice, the use of infinite sets or measures, or the use of the Fourier transform or inverse theorems from additive combinatorics.  It also gives explicit (but extremely poor) quantitative bounds.
\end{abstract}

\maketitle

\section{Introduction}

A famous theorem of van der Waerden \cite{vdw} in 1927 states the following.  

\begin{theorem}[Van der Waerden's theorem]\label{vwt}\cite{vdw}  For any integers $k,m \geq 1$ there exists an integer $N = N_{\vdW}(k,m) \geq 1$ such that every colouring $\c: \{1,\ldots,N\} \to \{1,\ldots,m\}$ of $\{1,\ldots,N\}$ into $m$ colours contains at least one monochromatic arithmetic progression of length $k$ (i.e. a progression in $\{1,\ldots,N\}$ of cardinality $k$ on which $\c$ is constant).
\end{theorem}

This theorem has by now several proofs; see \cite{walters} for a recent exposition of the original proof, as well as a proof of certain extensions of this theorem; for sake of completeness we present the original argument in an Appendix
(\S \ref{appendix-vwt}) below.  Another rather different proof can be found in \cite{shelah}.  This theorem was then generalized substantially in 1975 by Szemeredi \cite{szemeredi} (building upon earlier work in \cite{roth}, \cite{szemeredi-4}), answering a question of Erd\"os and Turan \cite{erdos}, as follows:

\begin{theorem}[Szemer\'edi's theorem]\label{szt}  For any integer $k \geq 1$ and real number $0 < \delta \leq 1$,
there exists an integer $N_{\SZ}(k,\delta) \geq 1$ such that for every $N \geq N_{\SZ}(k,\delta)$, every
set $A \subset \{1,\ldots,N\}$ of cardinality $|A| \geq \delta N$ contains at least one arithmetic progression of length $k$.
\end{theorem}

It is easy to deduce Van der Waerden's theorem from Szemer\'edi's theorem (with $N_{\vdW}(k,m) := N_{\SZ}(k,\frac{1}{m})$)
by means of the pigeonhole principle.  The converse implication however, is substantially less trivial.

There are many proofs already known for Szemer\'edi's theorem, which we discuss below; the main purpose of this paper is present yet another such proof.  This may seem somewhat redundant, but we will explain our motivation for providing another proof later in this introduction.

Remarkably, while Szemer\'edi's theorem appears to be solely concerned with arithmetic combinatorics, it has spurred
much further research in other areas such as graph theory, ergodic theory, Fourier analysis, and number theory; for instance it was a key ingredient in the recent result \cite{gt-primes} that the primes contain arbitrarily long
arithmetic progressions.  Despite the variety of proofs now available for this theorem, however, it is still regarded as
a very difficult result, except when $k$ is small.  The cases $k=1,2$ are trivial, and the case $k=3$ is by now relatively well understood (see \cite{roth}, \cite{furst-book}, \cite{rs}, \cite{soly}, \cite{bourg2}, \cite{heath}, \cite{bourgain-triples} for a variety of proofs).  The case $k=4$ also has a number of fairly straightforward proofs (see \cite{szemeredi-4}, \cite{roth-4}, \cite{gowers-4}, \cite{frankl}), although already the arguments here are more sophisticated than for
the $k=3$ case.  However for the case of higher $k$, only four types of proofs are currently known, all of which are rather deep.  The original proof of Szemer\'edi \cite{szemeredi} is highly combinatorial, relying on van der Waerden's theorem (Theorem \ref{vwt}) and the famous Szemer\'edi regularity lemma (which itself has found many other applications, see \cite{komlos} for a survey); it does provide an upper bound on $N_\SZ(k,\delta)$ but it is rather poor (of Ackermann type), due mainly to the reliance on the van der Waerden theorem and the regularity lemma, both of which have notoriously bad dependence of the constants.  Shortly afterwards, Furstenberg \cite{furst} (see also \cite{furstenberg}, \cite{furst-book}) introduced what appeared to be a completely different argument, transferring the problem into one of recurrence in ergodic theory, and solving that problem
by a number of ergodic theory techniques, notably the introduction of a Furstenberg tower (which is the analogue of the
regularity lemma).  This ergodic theory argument is the shortest and most flexible of all the known proofs, and has been the most successful at leading to further generalizations of Szemer\'edi's theorem (see for instance \cite{bergelson-leibman}, \cite{bl-hj}, \cite{FK1}, \cite{FK2}, \cite{FK3}).  On the other hand it uses the axiom of choice and as such does not provide any effective bounds for the quantity $N_\SZ(k,\delta)$.  The third proof is more recent, and is due to Gowers \cite{gowers} (extending earlier arguments in \cite{roth}, \cite{gowers-4} for small $k$).  It is based on combinatorics, Fourier analysis, and inverse arithmetic combinatorics (in particular multilinear versions of Freiman's theorem and the Balog-Szem\'eredi theorem).  It gives far better bounds on $N_\SZ(k,\delta)$
(essentially of double exponential growth in $\delta$ rather than Ackermann or iterated tower growth), but also requires
far more analytic machinery and quantitative estimates.  Finally, very recent arguments of Gowers \cite{gowers-reg} and R\"odl, Skokan, Nagle, Tengan, Tokushige, and Schacht \cite{rodl}, \cite{rodl2}, \cite{rstt}, \cite{nrs}, relying primarily on a hypergraph version of the Szemer\'edi regularity lemma, have been discovered; these arguments are somewhat similar in spirit to Szemer\'edi's original proof (as well as the proofs in \cite{rs}, \cite{soly} in the $k=3$ case
and \cite{frankl} in the $k=4$ case) but is conceptually somewhat more straightforward (once one accepts the need to work with hypergraphs instead of graphs, which does unfortunately introduce a number of additional technicalities).
Also these arguments can handle certain higher dimensional extensions of Szemer\'edi's theorem first obtained by 
ergodic theory methods in \cite{FK1}.

As the above discussion shows, the known proofs of Szemer\'edi's theorem are extremely diverse.  However, they
do share a number of common themes, principal among which is the establishment of a \emph{dichotomy between randomness and structure}.  Indeed, in an extremely abstract and heuristic sense, one can describe all the known proofs of Szemer\'edi's theorem collectively as follows.  Start with the set $A$ (or some other object which is a proxy for $A$, e.g. a graph,
a hypergraph, or a measure-preserving system).  For the object under consideration, define some concept of \emph{randomness} (e.g. $\eps$-regularity, uniformity, small Fourier coefficients, or weak mixing), and some concept of \emph{structure} (e.g. a nested sequence of arithmetically structured sets such as progressions or Bohr sets, or a partition of a vertex set into a controlled number of pieces, a collection of large Fourier coefficients, a sequence of almost periodic functions, a tower of compact extensions of the trivial $\sigma$-algebra, or a $k-2$-step nilfactor).  Obtain some sort of \emph{structure theorem} that splits the 
object into a structured component, plus an error which is random relative to that structured component.  To prove
Szemer\'edi's theorem (or a variant thereof), one then needs to obtain some sort of \emph{generalized von Neumann 
theorem} to eliminate the random error, and then some sort of \emph{structured recurrence theorem} for the structured component.  

Obviously there is a great deal of flexibility in executing the above abstract scheme, and this explains the large number
of variations between the known proofs of Szemer\'edi type theorems.  Also, each of the known proofs finds some parts of
the above scheme more difficult than others.  For instance, Furstenberg's ergodic theory argument requires some effort
(and the axiom of choice) to set up the appropriate proxy for $A$, namely a measure-preserving probability system, and 
the structured recurrence theorem (which is in this case a recurrence theorem for a tower of compact extensions) is 
also somewhat technical.  In the 
Fourier-analytic arguments of Roth and Gowers, the structured component is simply a nested sequence of long arithmetic progressions, which makes the relevant recurrence theorem a triviality; instead, almost all the difficulty resides in 
the structure theorem, or more precisely in enforcing the assertion that lack of uniformity implies a density increment 
on a smaller progression.  Gowers' more recent hypergraph argument is more balanced, with no particular step being
exceptionally more difficult than any other, although the fact that hypergraphs are involved does induce a certain
level of notational and technical complexity throughout.  Finally, Szemer\'edi's original argument contains significant
portions (notably the use of the Szemer\'edi regularity lemma, and the use of density increments) which fit very nicely into the above scheme, but also contains some additional combinatorial arguments to connect the various steps of the proof together.

In this paper we present a new proof of Szemer\'edi's theorem (Theorem \ref{szt}) which implements the above scheme in
a reasonably elementary and straightforward manner.  This new proof can best be described as a ``finitary'' or ``quantitative'' version of the ergodic theory proofs of Furstenberg \cite{furst}, \cite{furstenberg}, in which one stays entirely in the realm of finite sets (as opposed to passing to an infinite limit in the ergodic theory setting).  
As such, the axiom of choice is not used, and an explicit bound for $N_\SZ(k,\delta)$ is in principle possible\footnote{It may also be possible in principle to extract some bound for $N_\SZ(k,\delta)$ directly from the original Furstenberg argument via proof theory, using such tools as Herbrand's theorem; see for instance \cite{girard} where a similar idea is applied to the Furstenberg-Weiss proof of van der Waerden's theorem to extract Ackermann-type bounds from what is apparently a nonquantitative argument.  However, to the author's knowledge this program has not been carried out previously in the literature for the ergodic theory proof of Szemer\'edi proof.  Also we incorporate some other arguments in order to simplify the proof and highlight some new concepts (such as a new Banach algebra of uniformly almost periodic functions).}
(although the bound is extremely poor, perhaps even worse than Ackermann growth, and certainly worse than the bounds obtained by Gowers \cite{gowers}).  We also borrow some tricks and concepts
from other proofs; in particular from the proof of the Szemer\'edi regularity lemma we borrow the $L^2$ incrementation
trick in order to obtain a structure theorem with effective bounds, while from the arguments of Gowers \cite{gowers} we borrow the Gowers uniformity norms $U^{k-1}$ to quantify the concept of randomness.  One of our main innovations is to complement these norms with the (partially dual) \emph{uniform almost periodicity norms} $UAP^{k-2}$ to quantify the concept of an \emph{uniformly almost periodic function of order $k-2$}.  This concept will be defined rigorously later, but suffice to say for now that a model example of
a uniformly  almost periodic function of order $k-2$ is a finite polynomial-trigonometric sum $f: \Z_N \to \C$ of the form\footnote{Actually, these functions are a somewhat special class of uniformly almost periodic functions of order $k-2$, which one might dub the \emph{quasiperiodic functions of order $k-2$}.  The relationship between the two seems very closely related to the distinction in ergodic theory between $k-2$-step nilsystems and systems which contain polynomial eigenfunctions of order $k-2$; see \cite{furst-weiss}, \cite{host-kra2} for further discussion of this issue.  It is also closely related to the rather vaguely defined issue of distinguishing ``almost polynomial'' or ``almost multilinear'' functions from ``genuinely polynomial'' or ``genuinely multilinear'' functions, a theme which recurs in the work of Gowers \cite{gowers-4}, \cite{gowers}, and also in the theorems of Freiman and Balog-Szemer\'edi from inverse additive combinatorics which were used in Gowers' work.  It seems of interest to quantify and pursue these issues further.}
\begin{equation}\label{quasi}
 F(x) := \frac{1}{J} \sum_{j=1}^J c_j e( P_j(x) / N ) \hbox{ for all } x \in \Z_N,
\end{equation}
where $\Z_N := \Z/N\Z$ is the cyclic group of order $N$, $J \geq 1$ is an integer, the $c_j$ are complex numbers bounded in magnitude by 1, $e(x) := e^{2\pi i x}$, and the $P_j$ are polynomials of degree at most $k-2$ and with coefficients in $\Z_N$.  The uniform almost periodicity norms serve to quantify how closely a function behaves like \eqref{quasi}, 
and enjoy a number of pleasant properties, most notably that they form a Banach algebra; indeed one can think of these
norms as a higher order variant of the classical Wiener algebra of functions with absolutely convergent Fourier series.

The argument is essentially self-contained, aside from some basic facts such as the Weierstrass approximation theorem; the main external ingredient needed is van der Waerden's theorem (to obtain the recurrence theorem for uniformly almost periodic functions), and we supply the standard short proof of van der Waerden's theorem in Appendix \S \ref{appendix-vwt}.  As such, we do not
require any familiarity with any of the other proofs of Szemer\'edi's theorem, although we will of course discuss the relationship between this proof and the other proofs extensively in our remarks.
In particular we do not use the Fourier transform, or theorems from inverse arithmetic combinatorics such as Freiman's theorem or the Balog-Szemer\'edi theorem, and we do not explicitly use the Szemer\'edi regularity lemma either for graphs or hypergraphs (although the proof of that lemma has some parallels with certain parts of our argument here).  Also, while we do use the language of ergodic, measure, and probability theory,
in particular using the concept of conditional expectation with respect to a $\sigma$-algebra, we do so entirely
in the context of finite sets such as $\Z_N$; as such, a $\sigma$-algebra is nothing more than a finite partition of
$\Z_N$ into ``atoms'', and conditional expectation is merely the act of averaging a function on each atom\footnote{Readers familiar with the Szemer\'edi regularity lemma may see parallels here with the proof of that lemma.  Indeed one can
phrase the proof of this lemma in terms of conditional expectation; see \cite{tao:regularity}.}.  As such, we do not need
such results from measure theory as the construction of product measure (or conditional product measure, via Rohlin's lemma
\cite{rohlin}), which plays an important part of the ergodic theory proof, notably in obtaining the structure and recurrence theorems.  Also, we do not use the compactness of
Hilbert-Schmidt or Volterra integral operators directly (which is another key ingredient in Furstenberg's structure theorem), although we will still need a quantitative finite-dimensional version of this fact (see Lemmas \ref{tbp}, \ref{fra} below).
Because of this, our argument could technically be called ``elementary''.  However we will need a certain amount of
structural notation (of a somewhat combinatorial nature) in order to compensate for the lack of an existing body of
notation such as is provided by the language of ergodic theory.

In writing this paper we encountered a certain trade-off between keeping the paper brief, and keeping the paper well-motivated.  We have opted primarily for the latter; if one chose to strip away all the motivation and redundant arguments
from this paper one could in fact present a fairly brief proof of Theorem \ref{szt} (roughly 20 pages in length); see \cite{tao:abridged}.  We also had a similar trade-off between
keeping the arguments simple, and attempting to optimize the growth of constants for $N_{\SZ}(k,\delta)$ (which by the arguments here could be as bad as double-Ackermann or even triple-Ackermann growth); since it seems clear that the arguments here have no chance whatsoever to be competitive with the bounds obtained by Gowers' Fourier-analytic proof \cite{gowers} we have opted strongly in favour of the former.

\begin{remark}
Because our argument uses similar ingredients to the ergodic theory arguments, but in a quantitative finitary
setting, it seems likely that one could modify these arguments relatively easily to obtain quantitative finitary
versions of other ergodic theory recurrence results in the literature, such as those in \cite{FK1}, \cite{FK2}, \cite{FK3},
\cite{bergelson-leibman}, \cite{bl-hj}.  In many of these cases, the ordinary van der Waerden theorem
would have to be replaced by a more general result, but fortunately such generalizations are known to exist (see
e.g. \cite{walters} for further discussion).  In principle, the quantitative ergodic approach could in fact have a greater reach than the traditional ergodic approach to these problems; for instance, the recent establishment in \cite{gt-primes} that the primes contained arbitrarily long arithmetic progressions relied heavily on this quantitative ergodic point of view, and does not seem at this point to have a proof by traditional ergodic methods (or indeed by any of the other methods available for proving Szemer\'edi's theorem, although the recent hypergraph approach of Gowers \cite{gowers-reg} 
and of Rodl and Skokan \cite{rodl}, \cite{rodl2} seems to have a decent chance of being ``relativized'' to pseudorandom sets such as the ``almost primes'').  Indeed, some of the work used to develop this paper became incorporated into \cite{gt-primes}, and conversely some of the progress developed in \cite{gt-primes} was needed to conclude this paper.
\end{remark}

\begin{remark} It is certainly possible to avoid using van der Waerden's theorem explicitly in our arguments, for instance by incorporating arguments similar to those used in the proof of this theorem into the main argument\footnote{This is to some extent done for instance in Furstenberg's original proof \cite{furst}, \cite{furstenberg}.  A key component of that proof was to show that the multiple recurrence property was preserved under compact extensions.  Although it is not made explicit in those papers, the argument proceeds by ``colouring'' elements of the extension on each fiber, and using ``colour focusing'' arguments closely related to those used to prove van der Waerden's theorem (see e.g. \cite{walters}).  The relevance of van der Waerden's theorem and its generalizations in the ergodic theory approach is made more explicit in later papers, see e.g. \cite{furst-weiss}, \cite{bergelson-leibman}, \cite{bl-hj}, and also the discussion in \cite{walters}}.  A decreased reliance on van der Waerden's theorem would almost certainly lead to better bounds for $N_{\SZ}(k,\delta)$, for instance the Fourier-analytic arguments of Gowers \cite{gowers-4}, \cite{gowers} avoids this theorem completely and obtains bounds for $N_{\SZ}(k,\delta)$ which are far better than that obtained by any other argument, including ours.  However this would introduce additional arguments into our proof which more properly belong to the Ramsey-theoretic circle of ideas surrounding van der Waerden's theorem, and so we have elected to proceed by the simpler and ``purer'' route of using van der Waerden's theorem directly.  Also, as remarked above, the argument as presented here seems more able to extend to other recurrence problems.
\end{remark}

\begin{remark} Our proof of Szemer\'edi's theorem here is similar in spirit to the proof of the transference principle
developed in \cite{gt-primes} by Ben Green and the author which allowed one to deduce a Szemer\'edi theorem relative
to a pseudorandom measure from the usual formulation of Szemer\'edi's theorem; this transference principle also follows the same basic scheme used to prove Szemer\'edi's theorem (with Szemer\'edi's theorem itself taking on the role of the structured recurrence theorem).  Indeed, the two arguments were developed concurrently (and both were inspired, not only by each other, but by all four of the existing proofs of Szemer\'edi's theorem in the literature, as well as arguments from the much better understood $k=3,4$ cases); it may also be able to combine the two 
to give a more direct proof of Szemer\'edi's theorem relative to a pseudorandom measure.  There are two main differences however between our arguments here and those in \cite{gt-primes}.  Firstly, in the arguments here no pseudorandom measure is present.  Secondly, the role of structure in \cite{gt-primes} was played by the \emph{anti-uniform functions},
or more precisely a tower of $\sigma$-algebras constructed out of basic anti-uniform functions.  Our approach
uses the same concept, but goes further by analyzing the basic anti-uniform functions more carefully, and in fact
concluding that such functions are uniformly almost periodic\footnote{In \cite{gt-primes} the only facts required concerning these basic anti-uniform functions were that they were bounded, and that pseudorandom 
measures were uniformly distributed with respect to any sigma algebra generated by such functions.  This was basically because the argument in \cite{gt-primes} invoked Szemer\'edi's theorem as a ``black box'' to deal with this anti-uniform component, whereas clearly this is not an option for our current argument.} of a certain order $k-2$.
\end{remark}

\emph{Acknowledgements.}
This work would not have been possible without the benefit of many discussions with Hillel Furstenberg, Ben Green, Timothy Gowers, Bryna Kra, and Roman Sasyk for for explaining both the techniques and the intuition behind the various proofs of Szemer\'edi's theorem and related results in the literature, and for drawing the author's attention to various simplifications in these arguments.  Many of the ideas here were also developed during the author's collaboration with Ben Green, and we are particularly indebted to him for his suggestion of using conditional expectations and an energy increment argument to prove quantitative Szemer\'edi-type theorems.  We also thank Van Vu for much encouragement throughout this project, and Mathias Schacht for some help with the references.  The author also thanks Australian National University and Edinburgh University for their hospitality where this work was conducted.  The author is a Clay Prize Fellow and is supported by a grant from the Packard Foundation.

\section{The finite cyclic group setting}

We now begin our new proof of Theorem \ref{szt}.  Following the abstract scheme outlined 
in the introduction, we should begin by specifying what objects we shall use as proxies for the set $A$.
The answer shall be that we shall use non-negative bounded functions $f: \Z_N \to \R^+$ on 
a cyclic group $\Z_N := \Z/N\Z$.  In this section we set out some basic notation for such functions,
and reduce Theorem \ref{szt} to proving a certain quantitative recurrence property for these functions.

\begin{remark}  The above choice of object of study fits well with
the Fourier-based proofs of Szemer\'edi's theorem in \cite{roth}, \cite{roth-4}, \cite{gowers-4}, \cite{gowers}, at least
for the initial stages of the argument.  However in those arguments one eventually passes from $\Z_N$ to a smaller
cyclic group $\Z_{N'}$ for which one has located a density increment, iterating this process until randomness
has been obtained (or the density becomes so high that finding arithmetic progressions becomes
very easy).  In contrast, we shall keep $N$ fixed and use the group $\Z_N$ throughout the argument; it will be a certain
family of $\sigma$-algebras which changes instead.  This parallels the ergodic theory argument \cite{furst}, \cite{furstenberg}, \cite{furst-book}, but also certain variants of the Fourier argument such as \cite{bourg2}, \cite{bourgain-triples}.  It also fits well with the philosophy of proof of the Szemer\'edi regularity lemma.
\end{remark}

We now set up some notation.  We fix a large prime number $N$, and fix $\Z_N := \Z/N\Z$ to be the cyclic group of order $N$.  We will assume that $N$ is extremely large; basically, it will be larger than any quantity depending on any
of the other parameters which appear in the proof.  We will write $O(X)$ for a quantity bounded by $CX$
where $C$ is independent of $N$; if $C$ depends on some other parameters (e.g. $k$ and $\delta$), 
we shall subscript the $O(X)$ notation accordingly to indicate the dependence.  Generally speaking we will order
these subscripts so that the extremely large or extremely small parameters are at the right.

\begin{definition}
If $f: X \to \C$ is a function\footnote{Strictly speaking, we could give the entire proof of Theorem \ref{szt} using only real-valued functions rather than complex-valued, as is done in the ergodic theory proofs, thus making the proof slightly more elementary and also allowing for some minor simplifications in the notation and arguments.  However, allowing the functions to be complex valued allows us to draw more parallels with Fourier analysis, and in particular to discuss such interesting examples of functions as \eqref{quasi}.}, and $A$ is a finite non-empty subset of $X$,
we define the \emph{expectation of $f$ conditioning on $A$}\footnote{We have deliberately made this notation to coincide with the usual notations of probability $P(\Omega)$ and expectation $\E(f)$ for random variables to emphasize the probabilistic nature of many of our arguments, and indeed we will also combine this notation with the probabilistic one
(and take advantage of the fact that both forms of expectation commute with each other).  Note that one can think of $\E(f(x)|x \in A)$ as the conditional expectation of $f(x)$, where $x$ is a random variable with the uniform distribution on $X$, conditioning on the event $x \in A$.}
\[ \E(f|A) = \E(f(x) | x \in A) := \frac{1}{|A|} \sum_{x \in A} f(x)\] 
where $|A|$ of course denotes the cardinality of $A$.  If in particular $f$ is an indicator function $f = {\bf 1}_\Omega$
for some $\Omega \subseteq X$, thus $f(x) = 1$ when $x \in \Omega$ and $f(x) = 0$ otherwise, we write
$\P(\Omega | A )$ for $\E( {\bf 1}_\Omega | A)$.  Similarly, if $P(x)$ is an event depending on $x$,
we write $\P( P | A )$ for $\E( {\bf 1}_P | A)$, where ${\bf 1}_{P(x)} = 1$ when $P(x)$ is true and ${\bf 1}_{P(x)} := 0$ otherwise.
\end{definition}

We also adopt the following ergodic theory notation: if $f: \Z_N \to \R$ is a function, we define the shifts $T^n f: \Z_N \to R$ for any $n \in \Z_N$ or $n \in \Z$ by
$$ T^n f(x) := f(x+n),$$
and similarly define $T^n \Omega$ for any $\Omega \subset \Z_N$ by $T^n \Omega := \Omega - n$, thus $T^n 1_\Omega = 1_{T^n \Omega}$.  Clearly these maps are algebra homomorphisms (thus $T^n(fg) = (T^n f)(T^n g)$ and $T^n(f+g) = T^n f + T^n g$),
preserve constant functions, and also preserve expectation (thus $\E(T^n f | \Z_N) = \E(f | \Z_N)$).  They also form a group, thus $T^{n+m} = T^n T^m$ and $T^0$ is the identity, and are unitary with respect to the usual inner product
$\langle f, g \rangle := \E( f \overline{g} )$.  We shall also rely frequently\footnote{Of course, since the space of functions on $\Z_N$ is finite-dimensional, all norms are equivalent up to factors depending on $N$.  However in line with our philosophy that we only wish to consider quantities which are bounded uniformly in $N$, we think of these norms as being genuinely distinct.} on the Banach algebra norm
$$ \| f \|_{L^\infty} := \sup_{x \in \Z_N} |f(x)|$$
and the Hilbert space structure
$$ \langle f, g \rangle := \E( f \overline{g}); \quad \| f \|_{L^2} := \langle f, f\rangle^{1/2} = \E( |f|^2)^{1/2};$$
later on we shall also introduce a number of other useful norms, in particular
the \emph{Gowers uniformity norms} $U^{k-1}$ and
the \emph{uniform almost periodicity norms} $UAP^{k-2}$.

To prove Theorem \ref{szt}, it will suffice to prove the following quantitative recurrence version of that theorem.

\begin{definition}  A function $f: \Z_N \to \C$ is said to be \emph{bounded} if we have $\| f\|_{L^\infty} \leq 1$.
\end{definition}

\begin{theorem}[Quantitative recurrence form of Szemer\'edi's theorem]\label{szt-recur}  For any integer $k \geq 1$, any large prime integer $N \geq 1$, 
any $0 < \delta \leq 1$, and any non-negative bounded function $f: \Z_N \to \R^+$ with
\begin{equation}\label{f-density}
 \E(f | \Z_N) \geq \delta
\end{equation}
we have
\begin{equation}\label{recurrence-targ} \E( \prod_{j=0}^{k-1} T^{jr} f(x) | x,r \in \Z_N) \geq c(k,\delta)
\end{equation}
for some $c(k,\delta) > 0$.
\end{theorem}

\begin{remark}  This is the form of Szemer\'edi's theorem required in \cite{gt-primes}.  This result
was then generalized in \cite{gt-primes} (introducing a small error $o_{k,\delta}(1)$) 
by replacing the hypothesis that $f$ was bounded by the more general hypothesis
that $f$ was pointwise dominated by a \emph{pseudorandom measure}.  This generalization was crucial to obtain
arbitrarily long progressions in the primes.  We will not seek such generalizations here, although we do remark
that the arguments in \cite{gt-primes} closely parallel to the ones here.
\end{remark}

We now show how the above theorem implies Theorem \ref{szt}.

\begin{proof}[Proof of Theorem \ref{szt} assuming Theorem \ref{szt-recur}]  Fix $k, \delta$.  Let $N \geq 1$ be large,
and suppose that $A \subset \{1,\ldots,N\}$ has cardinality $|A| \geq \delta N$.  By Bertrand's postulate, we can find a
large prime number $N'$ between $kN$ and $2kN$.
We embed $\{1,\ldots,N\}$ in $\Z_{N'}$ in the usual manner, and let $A'$ be the image of $A$ under this embedding.
Then we have $\E({\bf 1}_{A'} | \Z_N ) \geq \delta/2k$, and hence by \eqref{recurrence-targ}
$$ \E( \prod_{j=0}^{k-1} T^{jr} 1_{A'}(x) | x,r \in \Z_{N'}) \geq c(k,\delta/2k),$$
or equivalently
$$ | \{ (x,r) \in \Z_{N'}: x, x+r, \ldots, x+(k-1)r \in A' \} | \geq c(k,\delta/2k) (N')^2.$$
Since $N' \geq kN$ and $A' \subset \{1,\ldots,N\}$, we see that $1 \leq x \leq N$ and $-N \leq r \leq N$ in the
above set.  Also we may remove the $r=0$ component of this set since this contributes at most $N$
to the above sum.  If $N$ is large enough, the right-hand side is still positive, and this 
implies that $A$ contains a progression $x, x+r, \ldots, x+(k-1)r$, as desired.
\end{proof}

\begin{remark}  One can easily reverse this implication and deduce Theorem \ref{szt-recur} from Theorem \ref{szt};
the relevant argument was first worked out by Varnavides \cite{varnavides}.  In the ergodic theory proofs, Szemer\'edi's theorem is also stated in a form 
similar to \eqref{recurrence-targ}, but
with $\Z_N$ replaced by an arbitrary measure-preserving system (and $r$ averaged over some interval $\{1,\ldots,N\}$
going to infinity), and the left-hand side was then shown to have positive limit inferior, rather than being bounded
from below by some explicit $c(k,\delta)$.  However these changes are minor, and again it is easy to pass from one
statement to the other, at least with the aid of the axiom of choice (see \cite{furstenberg} for some further 
discussion on this issue).
\end{remark}

It remains to deduce Theorem \ref{szt-recur}.  This task shall occupy the remainder of the paper.  

\section{Overview of proof}

We shall begin by presenting the high-level proof of Theorem \ref{szt-recur}, implementing the abstract scheme
outlined in the introduction.  

One of the first tasks is to define measures of randomness and structure
in the function $f$.  We shall do this by means of two families of norms\footnote{Strictly speaking, the $U^0$ and $U^1$ norms are not actually norms, and the $UAP^0$ norm can be infinite when $f$ is non-constant.  However, these issues will be irrelevant for our proof, and in the most interesting case $k \geq 3$ there are no such degeneracies.}: 
the \emph{Gowers uniformity norms}
$$ \| f \|_{U^0} \leq \| f \|_{U^1} \leq \ldots \leq \|f\|_{U^{k-1}} \leq \ldots \leq \| f \|_{L^\infty}$$
introduced in \cite{gowers} (and studied further in \cite{host-kra2}, \cite{gt-primes}) 
and a new family of norms, the \emph{uniform almost periodicity norms}
$$ \| f \|_{UAP^0} \geq \| f \|_{UAP^1} \geq \ldots \geq \| f \|_{UAP^{k-2}} \geq \ldots \geq \| f \|_{L^\infty}$$
which turn out to be somewhat dual to the Gowers uniformity norms.  We shall mainly rely on the
$U^{k-1}$ and $UAP^{k-2}$ norms; the other norms in the family are required only for mathematical induction 
purposes.  We shall define the Gowers uniformity and uniform almost periodicity norms rigorously in Sections \ref{unif-sec} 
and \ref{ap-sec}
respectively.  For now, we shall simply give a very informal (and only partially accurate) heuristic: a function 
bounded in $UAP^{k-2}$ will typically look something like the polynomially quasiperiodic function \eqref{quasi} where all the polynomials have degree at most $k-2$, 
whereas a function small in $U^{k-1}$ is something like a function which is ``orthogonal'' to all such quasiperiodic functions \eqref{quasi}.

Next, we state the three main sub-theorems which we shall use to deduce Theorem \ref{szt-recur}.  
The first sub-theorem, which is rather standard (and the easiest of the three to prove), asserts that Gowers-uniform functions (i.e. functions with small $U^{k-1}$ norm) are negligible for the purposes of computing \eqref{recurrence-targ}; it will be proven in Section \ref{unif-sec}.  

\begin{theorem}[Generalized Von Neumann theorem]\label{gvnt}\cite{gowers}  Let $k \geq 2$, and let $\lambda_1, \ldots, \lambda_k$ be distinct elements of $\Z_N$.  Then for any bounded functions
$f_1, \ldots, f_k: \Z_N \to \C$ we have
$$ |\E( \prod_{j=0}^{k-1} T^{\lambda_j r} f_j(x) | x,r \in \Z_N)| \leq \min_{1 \leq j \leq k} \| f_j \|_{U^{k-1}}.$$
\end{theorem}

\begin{remark}  As indicated, this part of the argument is based on the arguments of Gowers \cite{gowers}; however it is
purely combinatorial, relying on the Cauchy-Schwarz inequality rather than on Fourier analytic techniques (which occupy other parts of the argument in \cite{gowers}).  Variants of this theorem go back at least as far as
Furstenberg \cite{furst}; see also \cite{gt-primes}, \cite{host-kra2} for some variants of this theorem.  We remark that
the linear shifts $\lambda_j r$ can be replaced by more general objects such as polynomial shifts, after replacing the $U^{k-1}$ norm by a higher Gowers uniformity norm; this is implicit for instance in \cite{bergelson-leibman}.
\end{remark}

The second sub-theorem is a special case of the main theorem, and addresses the complementary situation to Theorem \ref{gvnt}, where $f$ is now uniformly almost periodic instead of Gowers-uniform; it will be proven in Section \ref{central-sec}.

\begin{theorem}[Almost periodic functions are recurrent]\label{recurrence}  Let $d \geq 0$ and $k \geq 1$ be integers, and let $f_{U^\perp}, f_{UAP}$ be non-negative bounded functions such that we have the estimates
\begin{align}
 \| f_{U^\perp} - f_{UAP} \|_{L^2} &\leq \frac{\delta^2}{1024 k}\label{central-0}\\
 \E( f_{U^\perp} | \Z_N ) &\geq \delta \label{central-1}\\
 \| f_{UAP} \|_{UAP^d} &< M\label{central-2}
\end{align}
for some $0 < \delta, M < \infty$.
Then we have
\begin{equation}\label{cdm}
 \E( \prod_{j=0}^{k-1} T^{\mu jr} f_{U^\perp}(x) | x \in \Z_N; 0 \leq r \leq N_1) \geq c_0(d,k,\delta,M)
\end{equation}
for some $c_0(k,\delta,d,M) > 0$ and all $\mu \in \Z_N$ and $N_1 \geq 1$.
\end{theorem}

\begin{remark} This argument is a quantitative version of certain ergodic theory arguments by Furstenberg and 
later authors, and is the only place where the van der Waerden theorem
(Theorem \ref{vwt}) is required.  It is by far the hardest component of the argument.  
In principle, the argument gives explicit bounds for $c_0(d,k,\delta,M)$
but they rely (repeatedly) on Theorem \ref{vwt} and are thus quite weak.
As mentioned earlier, we need this theorem only when $d=k-2$, but 
allowing $d$ to be arbitrary is convenient for the purposes of proving this theorem by induction.  It is important
that the quantity $\frac{\delta^2}{1024 k}$ used in the right-hand side of \eqref{central-0} does not depend on $M$.
This significantly complicates
the task of proving this theorem when $M$ is large, of course, since the error between $f_{U^\perp}$ and $f_{UAP}$ 
may seem to dominate whatever gain one can obtain from \eqref{central-2}. Nevertheless, one can cope with such large
errors by means of the machinery of $\sigma$-algebras and conditional expectation.  This ability to tolerate reasonably
large $L^2$ errors in this recurrence result is also crucially exploited in the ``Zorn's lemma'' step in the ergodic theory arguments, in which one shows that the limit of a chain of extensions with the recurrence property is also recurrent.  The parameters $\mu, N_1$ are technical and are needed to facilitate the inductive argument
used to prove this Theorem; ultimately we shall take $\mu := 1$ and $N_1 := N-1$.
\end{remark}

Finally, we need a structure theorem, proven in Section \ref{structure-sec}, that splits an arbitrary function into a Gowers-uniform component and an uniformly almost periodic component (plus an error).

\begin{theorem}[Structure theorem]\label{structure}  Let $k \geq 3$, and let $f$ be a non-negative bounded
function obeying \eqref{f-density} for some $\delta > 0$.  Then we can find a positive number $M = O_{k,\delta}(1)$,
a bounded function $f_U$, and non-negative bounded functions $f_{U^\perp}$, $f_{UAP}$ such that we have
the splitting
$$ f = f_U + f_{U^\perp}$$
and the estimates \eqref{central-0}, \eqref{central-1}, \eqref{central-2} with $d := k-2$, 
as well as the uniformity estimate
\begin{equation}\label{u-est}
 \| f_U \|_{U^{k-1}} \leq 2^{-k} c_0(k-2,k,\delta,M)
\end{equation}
where $c_0(d,k,\delta,M)$ is the quantity in \eqref{cdm}.
\end{theorem}

\begin{remark}  The subscripts $U$ and $U^\perp$ stand for \emph{Gowers uniform} and \emph{Gowers anti-uniform} respectively.  Thus this theorem asserts that while a general function $f$ need not have any uniformity properties
whatsoever, it can be decomposed into pieces which are either uniform in the sense of Gowers, or are instead uniformly almost periodicity, or are simply small in $L^2$.  This theorem is something of a hybrid between the Furstenberg 
structure theorem \cite{furstenberg} and the
Szemer\'edi regularity lemma \cite{szem-reg}.  A similar structure theorem was a key 
component to \cite{gt-primes}.  One remarkable fact here is that we could replace the quantity on the right-hand side
of \eqref{u-est} by an arbitrary positive function of $k$, $\delta$, $M$, at the cost of worsening the upper bound
on $M$.  The fact that the error tolerance in \eqref{central-0} does not go to zero as $M \to \infty$ is crucial in order to obtain this insensitivity to the choice of right-hand side of \eqref{u-est}.
\end{remark}

\begin{remark} Each of the above three theorems have strong parallels in the genuinely ergodic theory setting.  For instance, the analogues of the $U^d$ norms in that setting were worked out by Host and Kra \cite{host-kra2}, where the analogue of Theorem \ref{gvnt} was also (essentially) proven.  The structure theorem seems to correspond to the recent
discovery by Ziegler \cite{ziegler} of a universal characteristic factor for Szemer\'edi-type recurrence properties,
but with the role of the almost periodic functions of order $k-2$ replaced by the notion of a $k-2$-step nilsystem.
The recurrence theorem is very similar in spirit to $k-2$ iterations of the basic fact, established in \cite{furstenberg}, that recurrence properties are preserved under compact extensions (although our proof is not based on that argument, but instead on later colouring arguments such as the one in \cite{bergelson-leibman}).  One can also extend the definition of the Banach algebra $UAP^d$ defined below to the ergodic theory setting.  It seems of interest to pursue these connections
further, and in particular to rigorously pin down the relationship between almost periodicity of order $k-2$ and $k-2$-step nilsystems.
\end{remark}

Assuming these three theorems, we can now quickly conclude Theorem \ref{szt-recur}.

\begin{proof}[Proof of Theorem \ref{szt-recur}]  Let $f, k, \delta$ be as in Theorem \ref{szt-recur}.  We may take
$k \geq 3$ since the cases $k=1,2$ are trivial.  Let $M$, $f_U$, $f_{U^\perp}$, $f_{UAP}$ be as in Theorem \ref{structure}.
We can then split the left-hand side of \eqref{recurrence-targ} as the sum of $2^k$ terms of the form
$\E( \prod_{j=0}^{k-1} T^{jr} f_j(x) | x,r \in \Z_N)$, where each of the functions $f_0, \ldots, f_{k-1}$ is equal
to either $f_U$ or $f_{U^\perp}$.  The term in which all the $f_j$ are equal to $f_{U^\perp}$ is at least
$c_0(k,\delta,k-2,M)$ by Theorem \ref{recurrence} (taking $\mu := 1$ and $N_1 := N-1$).  
The other $2^k-1$ terms have magnitude at most
$\| f_U \|_{U^{k-1}} \leq 2^{-k} c_0(k,\delta,k-2,M)$ thanks to Theorem \ref{gvnt}.  Adding all this together we see that
$$\E( \prod_{j=0}^{k-1} T^{jr} f(x) | x,r \in \Z_N) \geq 2^{-k} c_0(k-2,k,\delta,M).$$
Since $M = O_{k,\delta}(1)$, the claim \eqref{recurrence-targ} follows.
\end{proof}

It remains to define the $U^{k-1}$ and $UAP^{k-2}$ norms properly, and
prove Theorems \ref{gvnt}, \ref{recurrence}, \ref{structure}.  This shall occupy the remainder of the paper.

\section{Uniformity norms, and the generalized von Neumann theorem}\label{unif-sec}

In this section we define the Gowers uniformity norms $U^d$ properly, and then prove Theorem \ref{gvnt}.  
The motivation for these norms comes from the \emph{van der Corput lemma}, which is very simple in
the context of the cyclic group $\Z_N$:

\begin{lemma}[Van der Corput Lemma]\label{vdc} For any function $f \in \Z_N \to \C$, we have
$$ |\E(f | \Z_N)|^2 = \E( \E( \overline{f} T^h f | \Z_N ) | h \in \Z_N ).$$
\end{lemma}

\begin{proof}  Expanding both sides the identity becomes
$$ \E( \overline{f}(x) f(y) | x,y \in \Z_N ) = \E( \overline{f(x)} f(x+h) | x,h \in \Z_N )$$
and the claim follows by the substitution $y = x+h$.
\end{proof}

Motivated by this lemma, we define

\begin{definition}[Gowers uniformity norms]\cite{gowers}\label{unif-def} Let $f: \Z_N \to \C$ be a function.  We define the $d^{th}$ Gowers uniformity norm $\| f \|_{U^d}$ recursively by
\begin{equation}\label{u0-def}
 \| f \|_{U^0} := \E(f | \Z_N)
\end{equation}
and
\begin{equation}\label{ud-def}
 \| f \|_{U^d} := \E( \| \overline{f} T^h f \|_{U^{d-1}}^{2^{d-1}} | h \in \Z_N )^{1/2^d}
\end{equation}
for all $d \geq 1$.
\end{definition}

\begin{example}  From Lemma \ref{vdc}, \eqref{u0-def}, \eqref{ud-def} we obtain the explicit formula
\begin{equation}\label{u1-def}
\| f \|_{U^1} = |\E(f | \Z_N)|.
\end{equation}
In particular, the $U^1$ norm (and hence all higher norms) are 
always non-negative.  The $U^2$ norm can also be interpreted as the $l^4$ norm of the Fourier coefficients 
of $f$ via the identity
\begin{equation}\label{l4} \| f \|_{U^2} = (\sum_{\xi \in \Z_N} |\E( f(x) e(-x\xi/N) | x \in \Z_N)|^4)^{1/4},
\end{equation}
though we will not need this fact here.  The higher $U^d$ norms do not seem have any particularly useful 
Fourier-type representations, however by expanding \eqref{ud-def} out recursively one can write the $U^d$ norm
as a sum of $f$ over $d$-dimensional cubes (see \cite{gowers}, \cite{gt-primes}, \cite{host-kra2} for further 
discussion of this). 
\end{example}

\begin{remark}  The $U^0$ and $U^1$ norms are not, strictly speaking, norms; the latter is merely a semi-norm, and the former is not a norm at all.  However, the higher norms $U^d$, $d \geq 2$
are indeed norms (they are homogeneous, non-degenerate, and obey the triangle inequality), and are also related to a certain $2^d$-linear inner product; see \cite{gowers}, \cite{gt-primes}, or \cite{host-kra2} for a proof of these facts 
(which we will not need here), with the $d=2$ case following directly from inspection of \eqref{l4}.  Also
one can show the inequality $\| f \|_{U^d} \leq \| f \|_{U^{d+1}}$ for any $d \geq 0$.  Thus for $k \geq 2$,
we have a rather interesting nested sequence of Banach spaces $U^{k-1}$ of functions $f: \Z_N \to \C$, equipped with 
the $U^{k-1}$ norm; these Banach spaces and their duals $(U^{k-1})^*$ were explored to a limited extent in \cite{gt-primes}, and we shall continue their study later in this paper.  Functions which are small in $U^2$ norm are termed \emph{linearly uniform} or \emph{Gowers-uniform of order 1}, and thus have small Fourier coefficients by \eqref{l4}; functions small in $U^3$ norm are \emph{quadratically uniform} or \emph{Gowers-uniform of order 2}, and so forth.  The terminology here is partly explained by the next example; again, see \cite{gowers}, \cite{gt-primes}, or \cite{host-kra2} for further discussion.
\end{remark}

\begin{example}  By induction\footnote{Actually, more is true: the $U^d$ norms of $f$ increase monotonically and converge to $\|f\|_{L^\infty}$ as $d \to \infty$, although the convergence can be quite slow and depends on $N$.  We will not prove this fact here.} we see that $\|f\|_{U^d} \leq \|f\|_{L^\infty}$ for all $d$; in particular we have $\|f\|_{U^d} \leq 1$ when $f$ is bounded.
We now present an example (which is, in fact, the only example up to scalar multiplication) in which
equality holds.  Let $P: \Z_N \to \Z_N$ be a polynomial with coefficients in $\Z_N$, and let $f(x) := e(P(x)/N)$.
Then one can show that $\| f \|_{U^d} = 1$ when $d \geq \deg(P)$, and $\|f\|_{U^d} = o_{\deg P}(1)$ when $d < \deg(P)$;
the former fact can be proven by induction and the trivial observation that for each fixed $h$, the polynomial $P(x+h) - P(x)$ has degree at most $\deg(P)-1$, while the latter fact also follows from induction, the above observation, and Lemma \ref{vdc}; we omit the details.  In fact one can improve the $o_{\deg P}(1)$ bound to $O_{\deg P}(N^{-1/2})$,
by using the famous Weil estimates.  By using the triangle inequality for $U^d$ (see e.g. \cite{gowers}, \cite{gt-primes}) one can also deduce similar statements for the polynomially quasiperiodic functions \eqref{quasi}.
\end{example}

One can easily verify by induction that the $U^d$ norms are invariant under shifts, thus $\|T^n f \|_{U^d} = \|f\|_{U^d}$,
and also invariant under dilations, thus if $\lambda \in \Z_N \backslash 0$ and $f_\lambda(x) := f(x/\lambda)$ then $\| f_\lambda \|_{U^d} = \| f \|_{U^d}$.

We can now prove the generalized von Neumann theorem.

\begin{proof}[Proof of Theorem \ref{gvnt}]
We induct on $k$.  When $k=2$ we use the fact that $(x,r) \mapsto (x + \lambda_1 r, x + \lambda_2 r)$ is
a bijection from $\Z_N^2$ to $\Z_N^2$ (recalling that $N$ is prime) to conclude that
$$ \E( \prod_{j=0}^{1} T^{\lambda_j r} f_j(x) | x,r \in \Z_N) = \E(f_1|\Z_N) \E(f_2|\Z_N)$$
and the claim then follows easily from \eqref{u1-def} and the boundedness of $f_1, f_2$.  
Now suppose that $k > 2$ and the claim
has already been proven for $k-1$.  By permuting the $\lambda_j$ if necessary we may assume that the minimum of 
the $\|f_j\|_{U^{k-1}}$ is attained when $j=0$.  By making the scaling $r \mapsto \lambda_0 r$ if
necessary we may assume that $\lambda_0 = 1$. By applying the expectation-preserving map
$T^{-\lambda_{k-1} r}$ (i.e. by subtracting $\lambda_{k-1}$ from each of the $\lambda_j$) we may assume that
$\lambda_{k-1}$ is zero.  The claim can now be written as
$$ |\E( f_{k-1}(x) \E( \prod_{j=0}^{k-2} T^{\lambda_j r} f_j(x) | r \in \Z_N) | x \in \Z_N)|
\leq \| f_0 \|_{U^{k-1}}.$$
By the Cauchy-Schwarz inequality and the boundedness of $f_{k-1}$, it suffices to prove that
$$ |\E( |\E( \prod_{j=0}^{k-2} T^{\lambda_j r} f_j(x) | r \in \Z_N)|^2 | x \in \Z_N)
\leq \| f_0 \|_{U^{k-1}}^2.$$
But from Lemma \ref{vdc} we have
\begin{align*}
\E( |\E( \prod_{j=0}^{k-2} T^{\lambda_j r} f_j(x) | r \in \Z_N)|^2 | x \in \Z_N)
&= \E( (\overline{\prod_{j=0}^{k-2} T^{\lambda_j r} f_j(x)}) (\prod_{j=0}^{k-2} T^{\lambda_j (r+h)} f_j(x)) | x,h,r \in \Z_N)\\
&= \E(\E( \prod_{j=0}^{k-2} T^{\lambda_j r} (\overline{f_j} T^{\lambda_j h})(x) | x,r \in \Z_N ) | h \in \Z_N ).
\end{align*}
On the other hand, from the induction hypothesis and the reduction to the case $\lambda_0 = 1$ we have
$$ |\E( \prod_{j=0}^{k-2} T^{\lambda_j r} (\overline{f_j} T^{\lambda_j h})(x) | x,r \in \Z_N )|
\leq \| f_0 T^h f_0 \|_{U^{k-2}}$$
for all $h \in \Z_N$.  Combining these two facts together, we obtain
$$ |\E( |\E( \prod_{j=0}^{k-2} T^{\lambda_j r} f_j(x) | r \in \Z_N)|^2 | x \in \Z_N)
\leq \E( \| f_0 T^h f_0 \|_{U^{k-2}} | h \in \Z_N )$$
and the claim follows from \eqref{ud-def} and H\"older's inequality.
\end{proof}

\begin{remark} The notion of Gowers uniformity considered here, namely that the $U^{k-1}$ norm is small, generalizes
the concept of \emph{pseudorandomness} or \emph{linear uniformity} in the $k=3$ case, which amounts to the assertion that all the Fourier coefficients of $f$ (except possibly for the zero coefficient) are small; this is the notion used for instance in \cite{roth}, \cite{bourg2}, \cite{heath}, \cite{bourgain-triples}.  This notion is essentially equivalent to the pair correlations of all the shifts $T^n f$ to be small on the average.  For higher $k$, this notion is insufficient
to obtain theorem such as Theorem \ref{gvnt}, see \cite{gowers-4}, \cite{gowers} for further discussion.  In Szemer\'edi's original arguments \cite{szemeredi-4}, \cite{szemeredi}, the appropriate concept of uniformity is provided by the notion of \emph{$\eps$-regularity}, which roughly corresponds to controlling all the $U^d$ norms for $d \leq C(\eps)$, while in the ergodic theory arguments of Furstenberg and later authors, the notion of uniformity used is that of \emph{weak mixing}, which roughly corresponds to controlling the $U^d$ norms for \emph{all} $d$.  Thus these notions of uniformity are significantly stronger than the one considered here, which fixes $d$ at $k-1$.  There is of course a cost to using such a strong notion of uniformity, and it is that one has to make the tower of structures extremely large in order to eventually attain such uniformity. In Szemer\'edi's regularity lemma, for instance, one is forced to lose constants which are of tower-exponential type in the regularity parameter $\eps$; see \cite{gowers-sz}.
In the ergodic theory arguments, the situation is even worse; the tower of invariant $\sigma$-algebras given by
Furstenberg's structure theorem (the ergodic theory analogue of Szemer\'edi's regularity lemma) can be as tall
as any countable ordinal, but no taller; see \cite{bf}.  
\end{remark}

\begin{remark} In the ergodic theory setting, one can also
define analogues of the $U^{k-1}$ norms, giving rise to the concept of invariant $\sigma$-algebras whose complement
consists entirely of functions which are \emph{Gowers-uniform of order $k-2$}; using this notion (which is much weaker than
weak mixing) it is possible to obtain a version of Furstenberg's structure theorem using only a tower of height $k-2$
(in particular, a tower of finite height).  Indeed, it was the author's discovery of this fact which led eventually
to the quantitative proof presented here; we have since learnt that this fact is essentially implicit in the work of
Host and Kra \cite{host-kra2} and Ziegler \cite{ziegler}.  
\end{remark}

\section{Almost periodic functions}\label{ap-sec}

Having defined the Gowers uniformity norms $U^{k-1}$ used for the generalized Von Neumann theorem, we now turn to defining
the dual concept of the uniform almost periodicity norms $UAP^{k-2}$ which we will need for both the recurrence
theorem and structure theorems.  Roughly speaking, if a function $F$ has a bounded $UAP^{k-2}$ norm, then it
should resemble a function of the form \eqref{quasi}, which we shall loosely refer to as a \emph{quasiperiodic function
of order $k-2$}.  To quantify this we make the following observation: if
$F$ is of the form \eqref{quasi}, then
\begin{equation}\label{shift-rep}
 T^n F = \E( c_{n,j} g_j | j \in J ),
\end{equation}
where $g_j$ is the bounded function $g_j(x) := e(P_j(x)/N)$, and $c_{n,j}(x)$ is the $c_{n,j}(x) = e( (P_j(x+n) - P_j(x)) / N)$.  The point here is that the dependence on $n$ on the right-hand side only arises through the functions $c_{n,j}$,
and those functions are of the form \eqref{quasi} but with degree $k-3$ instead of $k-2$.  Thus, the shifts of
an quasiperiodic functions of order $k-2$ can be written as linear combinations of fixed bounded functions $g_j$, where
the coefficients $c_{n,j}$ are not constant, but are instead quasiperiodic functions of one lower order\footnote{This observation was motivated by the use of \emph{relatively almost periodic functions} in
the ergodic theory arguments of Furstenberg \cite{furst}, \cite{furstenberg}, \cite{furst-book} and later authors.}.
We can pursue the same idea to define the $UAP^d$ norms recursively as follows.  

\begin{definition}[Banach algebras]  A space $A$ of functions on $\Z_N$, equipped with a norm $\|\|_A: A \to \R^+$, is said
to be a \emph{Banach algebra} if $A$ is a vector space, $\| \|_A$ is a norm (i.e. it is homogeneous, non-degenerate, 
and verifies the triangle invariant) which invariant under
conjugation $f \mapsto \overline{f}$, and  $A$ is closed under pointwise product with $\|fg\|_A \leq \|f\|_A \|g\|_A$
for all $f,g \in A$.  We also assume that
\begin{equation}\label{apd}
 \|F\|_{L^\infty} \leq \| F \|_{A}
\end{equation}
for all $F \in A$ (actually this property can be deduced from the pointwise product property and the finite-dimensionality of $A$).  We adopt the convention that $\|f\|_A = \infty$ if $f \not \in A$.  We say that $A$
is \emph{shift-invariant} if $\|f\|_A = \|T^n f \|_A$ for all $n$, and \emph{scale-invariant} if
$\|f_\lambda \|_A = \|f\|_A$ for all $\lambda \in \Z_N \backslash \{0\}$.
\end{definition}

\begin{definition}[Uniform almost periodicity norms]\label{ap-def}  If $A$ is a shift-invariant Banach algebra of functions on $\Z_N$, we define
the space $UAP[A]$ to be the space of all functions 
$F$ for which the orbit $\{ T^n F: n \in \Z \}$ has a representation of the form
\begin{equation}\label{representation}
 T^n F = M \E( c_{n,h} g_h ) \hbox{ for all } n \in \Z_N
\end{equation}
where $M \geq 0$, $H$ is a finite non-empty set, $g = (g_h)_{h \in H}$ is a collection of bounded functions, $c = (c_{n,h})_{n \in \Z_N, h \in H}$ is a collection of functions in $A$ with $\|c_{n,h}\|_{A} \leq 1$, 
and $h$ is a random variable taking values 
in $H$.  We define the norm $\|F\|_{UAP[A]}$ to be the infimum of $M$ over all possible representations
of this form.  
\end{definition}

\begin{remark}  Note that we are not imposing any size constraints on $H$, which could in fact get quite large (in fact one could allow $H$ to be infinite, i.e. $h$ could be a continuous random variable rather than a discrete one, without actually affecting this definition).  It turns out however that we will not need any information about $H$, or
more generally about the probability distribution of the random variable $h$.  The key point is that the 
Volterra operator $(c_h)_{h \in H} \mapsto \E( c_h g_h )$ will be a ``compact'' operator uniformly in choice of $h$ 
and $(g_h)_{h \in H}$.
\end{remark}

We first observe that the construction $A \mapsto UAP[A]$ maps shift-invariant Banach algebras to shift-invariant Banach algebras:

\begin{proposition}\label{ap-prop}  If $A$ is a shift-invariant Banach algebra, then so is $UAP[A]$.  Furthermore $UAP[A]$
contains $A$, and $\|f\|_{UAP[A]} \leq \|f\|_A$ for all $f \in A$.  Finally, if $A$ is scale-invariant then so is $UAP[A]$.
\end{proposition}

\begin{remark} This is a quantitative analogue of the well known fact in ergodic theory that the almost periodic functions form a shift-invariant algebra.
\end{remark}

\begin{proof}  It is easy to see that $UAP[A]$ is shift-invariant, conjugation-invariant, closed under scalar multiplication, preserves scale-invariance, and that the $UAP[A]$ norm is non-negative and homogeneous.  From \eqref{apd} and \eqref{representation} we see that $\|f\|_{L^\infty} \leq \|f\|_{UAP[A]}$ for all $f \in UAP[A]$, from which we deduce that the $UAP[A]$ norm is non-degenerate.  Also we easily verify that $UAP[A]$ cotnains $A$ with $\|f\|_{UAP[A]} \leq \|f\|_A$ for all $f \in A$.  Next, we show that $UAP[A]$ is closed 
under addition and that the $UAP^d$ norm enjoys the triangle inequality.  By homogeneity and nondegeneracy
it suffices to show that the unit ball is convex, i.e. if $F, F' \in UAP[A]$ are such that $\|F\|_{UAP[A]}, \| F'\|_{UAP[A]} < 1$ then $(1-\theta) F + \theta F' \in UAP[A]$ with $\| (1-\theta) F + \theta F' \|_{UAP[A]} < 1$ for all $0 \leq \theta \leq 1$.
By Definition \ref{ap-def} we can find non-empty finite sets $H, H'$, bounded functions $(g_h)_{h \in H}$ and $(g'_{h'})_{h' \in H'}$,
and functions $(c_{n,h})_{n \in \Z_N, h \in H}$ and $(c'_{n,h'})_{n \in \Z_N, h' \in H'}$ in $A$
and random variables $h$, $h'$ taking values in $H$ and $H'$ respectively such that we have the representations
\begin{equation}\label{reprep}
 T^n F = \E( c_{n,h} g_h ); \quad T^n F' = \E( c'_{n,h'} g'_{h'} | h' \in H' ) \hbox{ for all } n \in \Z_N
\end{equation}
and the estimates
$$ \| c_{n,h} \|_{A}, \| c'_{n,h'} \|_{A} \leq 1 \hbox{ for all }
n \in \Z_N, h \in H, h' \in H'.$$
Also, by relabeling $H'$ if necessary we may assume that $H$ and $H'$ are disjoint.  In such a case
we can concatenate $(c_{n,h})_{n \in \Z_N, h \in H}$ and $(c'_{n,h'})_{n \in \Z_N, h' \in H'}$ to a single 
function $(\tilde c_{n,\tilde h})_{n \in \Z_N, \tilde h \in H \cup H'}$ and similarly concatenate the $g$ functions to
$(\tilde g_{n,\tilde h})_{\tilde h \in H \cup H'}$.  If one then defines the random variable $\tilde h$ to equal $h$ with
probability $(1-\theta)$ and $h'$ with probability $\theta$ (or more precisely, the probability distribution of $\tilde h$
is $1-\theta$ times that of $h$ plus $\theta$ times that of $h'$ then one sees from linearity of expectation that
$$ T^n (F+F') = \E( \tilde c_{n,\tilde h} \tilde g_{\tilde h}  ) \hbox{ for all } n \in \Z_N$$
and the claim follows.

Next, we establish the algebra property.  By homogeneity and nondegeneracy again it suffices to show
that the unit ball is closed under multiplication.  To see this, start with \eqref{reprep}.  Without loss of generality
we may assume that the random variables $h, h'$ are independent (because it is only their individual distributions which
matter for \eqref{reprep}, not their joint distribution).  But in that case we have
$$ T^n (FF') = \E( c_{n,h,h'} g_{h,h'} ) \hbox{ for all } n \in \Z_N$$
where $c_{n,h,h'} := c_{n,h} c_{n,h'}$ and $g_{h,h'} := g_h g_{h'}$.  Since the product of two bounded functions
is a bounded function, the claim follows from the algebra property of $A$.
\end{proof}

Thanks to this proposition, we can define the $UAP^d$ norms recursively for $d \geq 0$, by setting
$UAP^0$ to be the trivial Banach algebra of all constant functions (equipped with the $L^\infty$ norm), and then setting
$UAP^d := UAP[UAP^{d-1}]$ for all $d \geq 1$.  Thus $UAP^d$ is a shift-invariant, scale-invariant Banach algebra
for all $d \geq 0$.

\begin{example}  If $F$ is of the form \eqref{quasi}, then one can
verify by \eqref{shift-rep} and induction that $F \in UAP^{k-2}$ with $\|F\|_{UAP^{k-2}} \leq 1$ (here $h$ is an element
of $\{1,\ldots,J\}$ chosen uniformly at random).  In particular, $UAP^1$ contains the Wiener algebra of functions with absolutely convergent Fourier series.  In our finitary setting of $\Z_N$, this implies that \emph{every} function lies in $UAP^1$ and hence in all higher $UAP^d$ norms, though the norm may grow with $N$ and thus be very large.  Note that the property of being in the Wiener algebra is substantially stronger than being \emph{almost periodic}, which is roughly equivalent to asking that the Fourier coefficients are summable in $l^{2-\eps}$ for some $\eps > 0$ rather than
being summable in $l^1$.  For further comparison, the property of being bounded in $(U^2)^*$, as discussed in \cite{gt-primes}, is stronger than being almost periodic but weaker than being bounded in $UAP^1$; it is equivalent to asking for the Fourier coefficients to be summable in $l^{4/3}$.  See \cite{gt-primes} for further discussion.
\end{example}

\begin{example} There are more subtle examples of almost periodic functions than the quasiperiodic ones.  One example is the function 
$$f(x) := e( [ax/N] b / N) \psi(ax/N \mod N) \psi(x/N \mod 1)$$ 
for some fixed $1 \leq a,b \leq N$, where 
$x$ is thought of as an integer from 1 to $N$, $[x]$ denotes the  integer part of $x$ and $\psi(x)$ is a smooth cutoff 
to the region $0.4 \leq x-[x] \leq 0.6$.  This function has an $UAP^1$ norm of $O(1)$ uniformly in $a, b, N$, but
the required representation of the form \eqref{representation} is not particularly obvious (for instance one can 
set $g_h$ to be various translations and modulations of $f$, and then $T^n f$ can be decomposed as an absolutely 
summable combination of the $g_h$ using smooth partitions of unity and Fourier series).  In this case, 
one can eventually work out that $f$ also has an absolutely convergent Fourier series; however things are even 
less clear for the function
\begin{equation}\label{fab}
f(x) := e( [ax/N] b x/ N) \psi(ax/N \mod N) \psi(x/N \mod 1), 
\end{equation}
which has an $UAP^2$ norm of $O(1)$ but seems to
have no particular resemblance with any quadratic phase function.  These ``generalized quadratic phase functions''
are related to 2-step nilsystems, which are known to not always admit quadratic eigenfunctions; see e.g. \cite{furst-weiss} for further discussion.  Intriguingly, hints of this ``generalized quadratic'' structure also emerge in the work of Gowers \cite{gowers-4}.  The situation here is still far from clear, though, and further study is needed. 
\end{example}

The structure theorem, Theorem \ref{structure}, can be viewed as some sort of duality relationship between
$UAP^{k-2}$ and $U^{k-1}$.  We now provide two demonstrations of this duality.   The first such demonstration
is rather simple, but is not actually used in the proof of Szemer\'edi's theorem; the second demonstration will
be to some extent a converse of the first and is one of the key components used to prove Theorem \ref{structure}.

\begin{proposition}[Uniformity is orthogonal to almost periodicity]\label{uap-dual}
Let $k \geq 2$.  For any functions $f, F$ with $F \in UAP^{k-2}$, we have
$$ |\langle f, F \rangle| \leq \|f\|_{U^{k-1}} \|F\|_{UAP^{k-2}}.$$
\end{proposition}

\begin{proof}  We induct on $k$.  When $k=2$ the claim follows from \eqref{u1-def} and the fact that $F$ is necessarily constant.  Now suppose that $k \geq 3$ and the claim has already been proven for $k-1$.  By homogeneity it suffices
to show that if $\|f\|_{U^{k-1}}, \|F\|_{UAP^{k-2}} < 1$ then $|\langle f, F \rangle| < 1$.

By Definition \ref{ap-def} we can find $H$, $g$, $c$, $h$ with the representation \eqref{representation}, with $g_h$
bounded and $\| c_{n,h} \|_{AP^{k-3}} \leq 1$ for all $n \in \Z_N$, $h \in H$.
Next, we use \eqref{representation} the unitary nature of $T^n$ to write 
$$\langle f,F \rangle = \langle T^n f, T^n F \rangle = \E( T^n f(x) \E( \overline{c_{n,h}}(x) \overline{g_h}(x)) | x \in \Z_N);$$
averaging over $n$ and rearranging we thus have
$$ \langle f, F \rangle = \E( \E( \E( T^n f(x) \overline{c_{n,h}}(x) | n \in \Z_N) \overline{g_h}(x) | x \in \Z_N )).$$
By the Cauchy-Schwarz inequality and the boundedness of the $g_h$, we thus have
$$ |\langle f, F \rangle| \leq \E(\E( |\E( T^n f(x) \overline{c_{n,h}}(x) | n \in \Z_N)|^2 | x \in \Z_N ))^{1/2}.$$
But from Lemma \ref{vdc} we have
$$ |\E( T^n f(x) \overline{c_{n,h}}(x) | n \in \Z_N)|^2
= \E( T^n ( \overline{f} T^r f )(x) c_{n,h}(x) \overline{c_{n+r,h}}(x) | n, r \in \Z_N)$$
whence
$$ |\langle f, F \rangle| \leq \E(\E( \langle \overline{f} T^r f, T^{-n}(\overline{c_{n,h}} c_{n+r,h}) \rangle | n,r \in \Z_N ))^{1/2}.$$
Since $UAP^{k-3}$ is a shift-invariant Banach algebra we have $\| T^{-n}(\overline{c_{n,h}} c_{n+r,h}) \|_{UAP^{k-3}} \leq 1$.
By the inductive hypothesis we thus have
$$ |\langle \overline{f} T^r f, T^{-n}(\overline{c_{n,h}} c_{n+r,h}) \rangle| \leq  
\| \overline{f} T^r f \|_{U^{k-2}},$$
whence
$$  |\langle f, F \rangle| \leq 
\E( \E( \| \overline{f} T^r f \|_{U^{k-2}} | n,r \in \Z_N )).$$
The outer expectation can be discarded since the quantity inside the expectation is deterministic.  We may similarly discard the redundant $n$ average.  Using Cauchy-Schwarz and \eqref{ud-def}, we thus obtain
$$ |\langle f, F \rangle| \leq \| f \|_{U^{k-1}} < 1$$
as desired.
\end{proof}

\begin{remark} One can use this Proposition to give an alternate proof of Theorem \ref{gvnt}, based on the
observation (easily verified by induction) that if $f_1, \ldots, f_{k-1}$ are bounded functions and
$\lambda_1, \ldots, \lambda_{k-1}$ are disjoint non-zero elements of $\Z_N$,
then $\E( \prod_{j=1}^{k-1} T^{\lambda_j r} f_j | r \in \Z_N )$ lies in $UAP^{k-2}$ with norm at most 1.
\end{remark}

\begin{remark} In the notation of \cite{gt-primes}, this shows that the $UAP^{k-2}$ norm is larger than or equal to the $(U^{k-1})^*$ norm.  However, the $UAP^{k-2}$ norm appears to be strictly stronger.  For instance, as observed in
\cite{gt-primes} the $(U^2)^*$ norm is the $l^{4/3}$ norm of the Fourier coefficients, and this norm does not form
a Banach algebra (in fact, it does not even control $L^\infty$ unless one loses a power of $N^{1/4}$), and so 
cannot be equivalent to the $UAP^1$ norm.
\end{remark}

We now give a partial converse to Proposition \ref{uap-dual}, which is the key to Theorem \ref{structure}.

\begin{lemma}[Lack of Gowers-uniformity implies correlation with a UAP function]\label{nonrandom}\cite{gt-primes}  
Let $f$ be a bounded function such that $\|f\|_{U^{k-1}} \geq \eps$ for some $k \geq 3$ and $\eps > 0$.
Then there exists a bounded function $F \in UAP^{k-2}$ with $\|F\|_{UAP^{k-2}} \leq 1$ such that
$|\langle f, F \rangle| \geq \eps^{2^{k-1}}$.
\end{lemma}

\begin{proof}  We need the concept of a \emph{dual function} from \cite{gt-primes};
the ergodic theory analogue of such functions have also been recently studied in \cite{host-kra2}, \cite{assani}.
For any function $f: \Z_N \to \C$ and any $d \geq 0$, we define the \emph{dual function
of order $d$} of $f$, denoted $\D_d(f)$, by the recursive formula
\begin{equation}\label{dual0-def} \D_0(f) := 1
\end{equation}
(i.e. $\D_0(f)$ is just the constant function 1) and
\begin{equation}\label{duald-def}
 \D_d(f) := \E( \overline{\D_{d-1}(\overline{f} T^h f)} T^h f | h \in \Z_N )
\end{equation}
for all $d \geq 1$.

We now claim the identity
$$ \langle f, \D_d(f) \rangle = \| f \|_{U^d}^{2^d}$$
for all $d \geq 0$.  When $d=0$ the claim follows from \eqref{dual0-def} and \eqref{u0-def}.  Now suppose inductively
that $d \geq 1$ and the claim has already been proven for $d-1$.  By \eqref{duald-def} (and the definition of the
inner product) we have
$$ \langle f, \D_d(f) \rangle = \E( \langle \D_{d-1}(\overline{f} T^h f), \overline{f} T^h f \rangle | h \in \Z_N ),$$
and the claim now follows from the inductive hypothesis and \eqref{ud-def}.

We thus set $F := \D_{k-1}(f)$.  It is clear from induction that $F$ is bounded; it remains to show that $F$
has $UAP^{k-2}$ norm less than 1.  Indeed, we make the more general claim that
$$ \| \D_d(f) \|_{UAP^{d-1}} \leq 1 \hbox{ for all bounded } f \hbox{ and all } d \geq 1.$$
When $d=1$ this is clear since $\D_d(f)$ is just the constant function $\E(f)$ in this case.  Now suppose
inductively that $d \geq 2$ and the claim has already been proven for $d-1$.  Applying $T^n$ to both sides of
\eqref{duald-def} and making the change of variables $h \leftarrow n+h$ we obtain
$$ T^n \D_d(f) := \E( c_{n,h} g_h | h \in \Z_N ), \hbox{ where } c_{n,h} :=
T^n \overline{\D_{d-1}(\overline{f} T^{h-n} f)} \hbox{ and } g_h := T^h f.$$
The functions $g_h$ are clearly bounded.  Since $UAP^{d-2}$ is a shift-invariant Banach algebra, we
see from inductive hypothesis that the
functions $c_{n,h}$ lie in $UAP^{d-2}$ with norm at most 1, and the claim follows (thinking of $h$ as an element of
$\Z_N$ chosen uniformly at random).
\end{proof}

\section{$\sigma$-algebras of almost periodic functions}

To prove both the recurrence theorem (Theorem \ref{recurrence}) or the structure theorem (Theorem \ref{structure}) 
it is convenient not just to work with almost periodic functions, but also
with certain \emph{$\sigma$-algebras} generated by them (in analogy with the ergodic theory arguments).
Thus in this section we develop the theory of such $\sigma$-algebras.

We begin by recalling the definition of a 
$\sigma$-algebra (not necessarily associated with an almost periodic function).  This very useful concept is of course
equivalent in the finitary setting of $\Z_N$ to the more familiar notion of a \emph{partition}
of $\Z_N$, but we will retain the language of probability theory in order to maintain the analogy with the ergodic
theory arguments, and in order to benefit from such useful concepts as conditional expectation, orthogonality, measurability, energy, and so forth.  See \cite{tao:regularity} for a further discussion of the connection 
between the Szemer\'edi regularity lemma, partitions, and conditional expectation with respect to $\sigma$-algebras.

\begin{definition}[$\sigma$-algebras] A \emph{$\sigma$-algebra} $\B$ in $\mathbb{Z}_N$ is any collection
of subsets of $\mathbb{Z}_N$ which contains the empty set $\emptyset$ and the full set $\Z_N$, and is closed under
complementation, unions and intersections.  We define the \emph{atoms} of a $\sigma$-algebra to be the minimal non-empty elements of $\B$ (with respect to set inclusion); it is clear that the atoms in $\B$ form a partition of $\mathbb{Z}_N$, and $\B$ consists precisely of arbitrary unions of its atoms (including the empty union $\emptyset$); thus there is a one-to-one correspondence between $\sigma$-algebras and partitions of $\Z_N$.
A function $f: \Z_N \to \C$ is said to be \emph{measurable} with respect to a $\sigma$-algebra $\B$ if all the level 
sets of $f$ lie in $\B$, or equivalently if $f$ is constant on each of the atoms of $\B$.  

We define $L^2(\B) \subseteq L^2(\mathbb{Z}_N)$ to be the closed subspace of the Hilbert space $L^2(\mathbb{Z}_N)$ consisting of $\B$-measurable functions.  We can then define the conditional expectation operator $f \mapsto \E(f|\B)$ to
be the orthogonal projection of $L^2(\mathbb{Z}_N)$ to $L^2(\B)$.  An equivalent definition of conditional expectation is
$$ \E(f|\B)(x) := \E( f(y) | y \in \B(x) )$$
for all $x \in \mathbb{Z}_N$, where $\B(x)$ is the unique atom in $\B$ which contains $x$.  It is clear that conditional
expectation is a linear self-adjoint orthogonal projection on $L^2(\mathbb{Z}_N)$, preserves non-negativity, expectation, and constant functions.  In particular it maps bounded functions to bounded functions.  If $\E(f|\B)$ is zero we say
that $f$ is \emph{orthogonal to} $\B$.

If $\B$, $\B'$ are two $\sigma$-algebras, we use $\B \vee \B'$ to denote the $\sigma$-algebra generated by $\B$ and $\B'$ (i.e. the $\sigma$-algebra whose atoms are the intersections of atoms in $\B$ with atoms in $\B'$).  
\end{definition}

Observe that when $\B$ is the trivial $\sigma$-algebra $\{\emptyset,\Z_N\}$ then the conditional expectation
$\E(f|\B)$ is just the constant function equal to the ordinary expectation $\E(f|\Z_N)$.  Every $\sigma$-algebra
induces a unique orthogonal decomposition $f = \E(f|\B) + (f-\E(f|\B))$ of a function $f$ into the component
$\E(f|\B)$ measurable with respect to $\B$, and the component $f - \E(f|\B)$ orthogonal to $\B$.  
More generally, if $\B$ is a subalgebra of $\B'$ (thus $\B'$ is \emph{finer} than $\B$, or $\B$ is \emph{coarser}
than $\B'$) then we can orthogonally decompose the finer expectation $\E(f|\B') = \E(f|\B) + (\E(f|\B') - \E(f|\B))$
into the coarser expectation, and a component measurable in $\B'$ but orthogonal to $\B$.

We now show that each almost periodic function generates a well-behaved $\sigma$-algebra at every scale $\eps$.

\begin{proposition}[UAP functions generate a compact $\sigma$-algebra]\label{ap-generate}  Let $d \geq 0$, let $G \in UAP^d$ be such that $\| G \|_{UAP^d} \leq M$ for some $M > 0$, and let $\eps > 0$.  Then there exists a $\sigma$-algebra
$\B_{\eps}(G) = \B_\eps(G,d)$ consisting of at most $O_{M,\eps}(1)$ atoms, such that we have the following two properties:

\begin{itemize}

\item ($G$ lies in its own $\sigma$-algebra) We have the approximation property
\begin{equation}\label{g-approx}
 \| G - \E( G | \B_\eps(G) ) \|_{L^\infty} = O(\eps).
\end{equation}
Similarly if $\B_\eps(G)$ is replaced by any finer $\sigma$-algebra.

\item (Approximation by almost periodic functions)  For any bounded non-negative function $f$ which is measurable
in $\B_\eps(G)$, and any $\delta > 0$, there exists a bounded non-negative function $f_{UAP} \in UAP^d$
such that
\begin{equation}\label{efg-decomp-1}
\| f - f_{UAP} \|_{L^2} \leq \delta
\end{equation}
and
\begin{equation}\label{fap-ap-1}
\| f_{UAP} \|_{UAP^{d}} = O_{M, \eps, \delta}(1).
\end{equation}

\end{itemize}
\end{proposition}

\begin{remark}  As the proof shows, the above Proposition in fact holds if $UAP^d$ is replaced by any other Banach 
algebra.
\end{remark}

\begin{proof}  We shall prove this by constructing $\B_{\eps}(G)$ using randomized level sets (or ``generalized Bohr sets'' of $G$), using some ideas from \cite{gt-primes}.   Let $S := \{ z \in \C: -1/2 \leq |\Re(z)|, |\Im(z)| < 1/2 \}$ be the unit square in the complex plane, and let $\Z[i] := \{ a+bi: a,b \in \Z\}$ denote the 
Gaussian integers.  Let $\alpha \in S$ be a complex number chosen uniformly at random from $S$.
We can then define the $\sigma$-algebra
$\B_{\eps,\alpha}(G)$ to be the algebra whose atoms are the sets 
$\{G^{-1}(\eps ( S + \zeta + \alpha)): \zeta \in \Z[i] \}$.
It will suffice to show that with positive probability this algebra $\B_{\eps,\alpha}(G)$ is a candidate 
for $\B_\eps(G)$ (with the bounds in \eqref{g-approx}, \eqref{fap-ap} uniform in $\alpha$).

The bound \eqref{g-approx} is clear since $G$ is constrained to lie in a square of diameter $O(\eps)$ on
each atom of $\B_\eps$ (and hence on each atom of any finer $\sigma$-algebra).  
Since $\|G\|_{L^\infty} \leq \|G\|_{UAP^d} \leq M$, $G$ takes values in a ball of radius $O(M)$, and thus the number of atoms
in $\B_{\eps}$ is indeed $O_{M,\eps}(1)$ as claimed.  Now we turn to the approximation property.  It will suffice to prove that for each $\delta > 0$
and every $\eta > 0$, the approximation property is true (with the bound in \eqref{fap-ap} allowed to depend on $\eta$) with probability at least $1 - \eta$, since one can then set $\delta := 2^{-n}$ and $\eta := \delta/2$ for $n=1,2,\ldots$ (for instance) and conclude the claim is true for all $\delta$ with positive probability.

Now fix $\delta, \eta$.  Since every bounded non-negative functions $f$ is a convex combination of indicator functions
of the form ${\bf 1}_\Omega$ where $\Omega \in \B_{\eps,\alpha}(G)$, and the number of such functions is $O_{M,\eps}(1)$
(since $\B_{\eps,\alpha}(G)$ only $O_{M,\eps}(1)$ atoms), it suffices (after shrinking $\eta$ appropriately) to prove
the claim for a single indicator function $f := {\bf 1}_\Omega$.

Fix $\Omega$; we can then write
$$ f = {\bf 1}_{\Omega} = {\bf 1}_W( \eps^{-1} G - \alpha)$$
where $W \subset \C$ is some union of $O_{X}(1)$ translates of the unit cube $S$.

Let $0 < \sigma \ll 1$ be a small number (it will eventually be much smaller than $\delta$, $\eta$, or $\eps$)
to be chosen later.  Let $\partial W_\sigma$ be the $\sigma$-neighbourhood of
the boundary $\partial W$ of $W$.  By Urysohn's lemma combined with the Weierstrass approximation theorem
(and the fact that $G = O(M)$) we can write
\begin{equation}\label{om-decomp}
 f = P( \eps^{-1} G - \alpha, \overline{\eps^{-1} G - \alpha} )
+ O( \sigma ) + O(1) {\bf 1}_{\partial W_{\sigma}}( \eps^{-1} G - \alpha)
\end{equation}
for some polynomial $P = P_{W, M, \eps, \sigma}$ of two complex variables.  
Denote the first term on the right-hand side of
\eqref{om-decomp} by
$f_{UAP}$, then from $UAP^d$ hypothesis on $G$ and the Banach algebra nature of $UAP^d$ we have
\begin{equation}\label{fap-ap-d}
 \| f_{UAP} \|_{UAP^d} = O_{M, P}(1) = O_{M, \eps, \sigma}(1),
\end{equation}
which will give \eqref{fap-ap-1} once $\sigma$ is selected properly at the end of the argument.

Now consider the third term on the right-hand side of \eqref{om-decomp}.  Observe that
\begin{align*}
\| {\bf 1}_{\partial W_{\sigma}}( \eps G - \alpha) \|_{L^2}^2
&= \P( \eps^{-1} G(x) - \alpha \in \partial W_\sigma | x \in \Z_N ) \\
&\leq \P( \eps^{-1} G(x) - \alpha \in \partial S_\sigma + \zeta \hbox{ for some } \zeta \in \Z[i]
| x \in \Z_N )
\end{align*}
where $\partial S_\sigma$ is the $\sigma$-neighbourhood of the boundary $\partial S$ of the unit square.  Observe
that as $\alpha$ varies over $S$, the event
$$ \eps^{-1} G(x) - \alpha \in \partial S_\sigma + \zeta \hbox{ for some } \zeta \in \Z[i]$$
has probability $O(\sigma)$ regardless of what $\eps^{-1} G(x)$ is.  We thus have
$$ \E(\| {\bf 1}_{\partial W_{\sigma}}( \eps G - \alpha) \|_{L^2}^2)
\leq O(\sigma).$$
By Markov's inequality, we thus see that the expression inside the expectation is $O_{\eta}(\sigma)$ with probability
at least $1-\eta$.  Inserting this into \eqref{om-decomp} we obtain
$$ \| f - f_{UAP} \|_{L^2} = O(\sigma) + O_{\eta}(\sigma^{1/2})$$
and the claim \eqref{efg-decomp-1} follows by setting $\sigma$ sufficiently small depending on $\delta$ (and
then \eqref{fap-ap-1} will hold by \eqref{fap-ap-d}).
\end{proof}

Henceforth we shall fix an assignment of a $\sigma$-algebra $\B_{\eps}(G) = \B_{\eps}(G,d)$ with the above properties
for each almost periodic function $G \in UAP^d$ and each $\eps > 0$.  Note that while we did use a randomization argument
here, it is possible to make such an assignment constructive, for instance by well-ordering all the $\sigma$-algebras of
$\Z_N$ in some constructive way and then choosing the minimal algebra which obeys the above properties (with the exact choice of bounds in \eqref{fap-ap}, etc. held fixed).  Thus we do not require the axiom of choice at this step (or indeed at any step in this argument).  For similar reasons we may ensure that this procedure is shift-invariant, in the sense that
\begin{equation}\label{shift-invariant}
T^n \B_\eps(G) = \B_\eps(T^n G) \hbox{ for all } n \in \Z_N,
\end{equation}
where $T^n \B := \{ T^n \Omega: \Omega \in \B \}$ is the $\sigma$-algebra $\B$ shifted backwards by $n$.  This shift
invariance amounts to making sure the same $\alpha$ is chosen for all the shifts $T^n G$ of a fixed function $G$,
which is easy enough to ensure since the constraints needed for $\alpha$ are independent of the choice of $n$.

The above Proposition pertained to a $\sigma$-algebra generated by a single almost periodic function, but we can
easily extend it to algebras generated by multiple functions as follows.

\begin{definition}[Compact $\sigma$-algebras]\label{complexity-def}  Let $d \geq 0$ and $X \geq 0$.  
A $\sigma$-algebra $\B$ is said to 
be \emph{compact of order $d$ and complexity at most $X$} if it has the form
\begin{equation}\label{b-def}
 \B = \B_{\eps_1}(G_1) \vee \ldots \vee \B_{\eps_K}(G_K)
\end{equation}
for some $0 \leq K \leq X$, some $\eps_1, \ldots, \eps_K \geq \frac{1}{X+1}$, and some $G_1, \ldots, G_K \in UAP^{d-1}$
with norm $\| G_j \|_{UAP^{d}} \leq X$ for all $1 \leq j \leq K$.  We define the \emph{$d$-complexity} (or simply
\emph{complexity}) of a $\sigma$-algebra to be the minimal $X$ for which one has the above representation, or
$\infty$ if no such represenation exists. In particular, the trivial $\sigma$-algebra
$\B = \{ \emptyset, \Z_N\}$ is compact of order $d$ with complexity 0.
\end{definition}

\begin{remark} 
The terminology is motivated here by ergodic theory, see e.g. \cite{furstenberg}; a compact $\sigma$-algebra of order
$d$ here corresponds in the ergodic setting, roughly speaking, to a 
tower of height $d$ of compact extensions of the trivial algebra.  The complexity $X$ is a rather artificial quantity
which we use as a proxy for keeping all the quantities used to define $\B$ under control.
\end{remark}

The key property we need concerning these $\sigma$-algebras is that the measurable functions of $\sigma$-algebras
which are compact of order $d$ are (with high probability) well approximated by almost periodic functions of
order $d$:

\begin{proposition}[UAP functions are dense in compact $\sigma$-algebras]\label{compact-sigma}  Let 
$d \geq 0$, $X \geq 0$, and let $\B$ be a $\sigma$-algebra which is compact of order $d$ and complexity at most $X$.
Let $f$ be a bounded non-negative function which is measurable with respect to $\B$, and let $\delta > 0$.
Then we can find a
bounded non-negative $f_{UAP} \in UAP^{d}$ such that
\begin{equation}\label{efg-decomp}
\| f - f_{UAP} \|_{L^2} \leq \delta
\end{equation}
and
\begin{equation}\label{fap-ap}
\| f_{UAP} \|_{UAP^{d}} = O_{d, \delta, X}(1).
\end{equation}
\end{proposition}

\begin{proof}  We first verify the claim when $f$ is the indicator function ${\bf 1}_A$ of an atom $A$ of $\B$.
From Definition \ref{complexity-def} we can expand $\B$ in the form \eqref{b-def}, and hence we can
write $A = A_1 \cap \ldots \cap A_K$ where each $A_j$ is an atom in $\B_{\eps_j}(G_j)$.  From Proposition \ref{ap-generate}
and the bounds on $\eps_j$, $G_j$, $K$ coming from Definition \ref{complexity-def},
we can find bounded non-negative functions $f_{UAP,j} \in UAP^d$ for all $1 \leq j \leq K$ such that
$$ \| {\bf 1}_{A_j} - f_{UAP,j} \|_{L^2} \leq \delta/K$$
and
$$ \| f_{UAP,j} \|_{L^2} = O_{d, \delta/K, \eps_j, X}(1) = O_{d, \delta, X}(1).$$
Since ${\bf 1}_{A_j}$ and $f_{UAP,j}$ are both bounded and non-negative we have the elementary pointwise inequality
$$ |\prod_{j=1}^K {\bf 1_{A_j}} - \prod_{j=1}^K f_{UAP,j}|
\leq \sum_{j=1}^K |{\bf 1_{A_j}} - f_{UAP,j}|$$
and hence if we set $f_{UAP} := \prod_{j=1}^K f_{UAP,j}$ then \eqref{efg-decomp} follows from the triangle inequality,
and \eqref{fap-ap} follows the Banach algebra nature of $UAP^d$.  Since $f_{UAP}$ is clearly bounded and non-negative,
the claim follows.

Now suppose $f$ is an arbitrary bounded non-negative function measurable with respect to $\B$.  Then we can
write $f = \sum_A c_A {\bf 1_A}$ where $A$ ranges over the atoms of $\B$ and $0 \leq c_A \leq 1$ are constants.
Let $\sigma = \sigma(d,\delta,X) > 0$ be a small number to be chosen later, then by the preceding discussion we can find bounded 
non-negative $f_{UAP,A} \in UAP^d$ for all $A$ such that
$$ \| {\bf 1_A} - f_{UAP,A} \|_{L^2} \leq \sigma$$
and
$$ \| f_{UAP,A} \|_{UAP^d} = O_{d,X,\sigma}(1).$$
If we then set $\tilde f_{UAP} := \sum_A c_A f_{UAP,A}$ and observe from Proposition \ref{ap-generate} that $\B$ contains at most $O_{d,X}(1)$ atoms, we thus have
$$ \| f - \tilde f_{UAP} \|_{L^2} \leq O_{d,X}(\sigma)$$
and
$$ \| \tilde f_{UAP} \|_{UAP^d} = O_{d,X,\sigma}(1).$$
We are however not done yet, because while $f_{UAP}$ is non-negative, it is not bounded by 1; instead we have
a bound of the form $0 \leq f_{UAP}(x) \leq O_{d,X}(1)$.  To fix this we need a real-valued 
polynomial $P(x) = P_{d,\delta,X}(x)$
such that
$$ |P(x) - \max(x,1)| \leq \delta/2 \hbox{ and } 0 \leq P(x) \leq 1 \hbox{ for all } x \hbox{ in the range of } f_{UAP};$$
such a polyomial exists by the Weierstrass approximation theorem.  If we then set $f_{UAP} := P(\tilde f_{UAP})$,
then $f_{UAP}$ is bounded and non-negative, and we have
$$ \| f_{UAP} - \max(\tilde f_{UAP},1) \|_{L^2} \leq \delta/2$$
and (since $UAP^d$ is a Banach algebra)
$$ \| f_{UAP} \|_{UAP^d} = O_{d,\delta,X,\sigma}(1).$$
On the other hand, since $f$ is bounded above by 1, we have
$$ \| f - \max(\tilde f_{UAP},1) \|_{L^2} \leq \| f - \tilde f_{UAP} \|_{L^2} \leq O_{d,X}(\sigma),$$
and the claims then follow from the triangle inequality if $\sigma$ is chosen sufficiently small depending
on $d$, $\delta$, $X$.
\end{proof}

\section{The energy incrementation argument}

The proof of the recurrence theorem (Theorem \ref{recurrence}) and the structure theorem (Theorem \ref{structure})
relies not only on $\sigma$-algebras of almost periodic functions, which we constructed in the previous section,
but also on the notion of the \emph{energy} of a $\sigma$-algebra
with respect to a collection of functions, and of the recursive \emph{energy incrementation argument} which 
we will need to prove both the recurrence theorem and the structure theorem.  This energy incrementation argument,
which was inspired by the proof of the Szemer\'edi regularity lemma (see e.g. \cite{szem-reg}), is perhaps
one of the most important aspects of this argument, but unfortunately is also the one which causes the Ackermann-type
(or worse) blowup of bounds.  It is the counterpart of the more well-known \emph{density incrementation argument}
which appears in several proofs of Szemer\'edi's theorem (starting with Roth's original argument \cite{roth}, but see also \cite{gowers-4}, \cite{gowers}, \cite{heath}, \cite{bourgain-triples}, \cite{szemeredi-4}, \cite{roth-4}, \cite{szemeredi}).  In that strategy one passes from the original set $\{1,\ldots,N\}$ to a decreasing sequence of
similarly structured subsets (e.g. arithmetic progressions or Bohr sets) while forcing the density $\delta$ of the set $A$ to increase as one progresses along the sequence; eventually one finds enough ``randomness'' to obtain an arithmetic progression.  The hope is to show this algorithm terminates successfully by using the trivial fact that the density is always bounded above by 1.  To do this, it is important that the density increment depend only on $\delta$, and 
not on other parameters such as $N$ or the complexity of the structured subset.  This rather stringent requirement
on the density increment is one cause of technical complexity and length in several of the arguments mentioned above.

In our situation, the role of ``structured subset'' will be played by a $\sigma$-algebra generated by almost periodic functions, and the role of density played by the energy of that $\sigma$-algebra.  This energy will automatically increase as the $\sigma$-algebra gets finer, and is also automatically bounded.  Once again, however, the energy increment may 
be very small, depending for instance on the complexity of the $\sigma$-algebra, and this algorithm may once again fail 
to terminate.  This problem also appears in the ergodic theory setting, in the context of an infinite tower of
$\sigma$-algebra extensions; to resolve this one must show that the supremum of any tower of extensions with
the recurrence property also has the recurrence property.  This appears difficult since the $\sigma$-algebras in 
this tower may become arbitrarily complex, and the lower bound obtained by the recurrence property may go to zero
as one approaches the supremum of the tower.  Nevertheless, one can conclude the argument, basically by observing
that any measurable function in the supremum of the tower can be approximated in $L^2$ norm (say) by a measurable
function in some finite component of this tower, and a simple argument then allows one to deduce recurrence for the former function from recurrence from the latter function regardless of how small the recurrence bound is for the latter;
see e.g. \cite{furstenberg} for an example of this.  This may serve to explain why we have the error tolerance
\eqref{central-0} in the recurrence theorem.

We begin by defining the \emph{energy} of a $\sigma$-algebra (relative to some fixed collection
of functions); this can be thought of as somewhat
analogous to the more standard notion of the \emph{entropy} of an algebra in both information theory and ergodic theory, but the energy will be adapted to a specific fixed collection of functions $f_1, \ldots, f_m$, whereas the entropy
is in some sense concerned with all possible functions at once.

\begin{definition}[Energy]\label{energy-def}
Given a $m$-tuple $f = (f_1, \ldots, f_m)$ of functions $f_j: \Z_N \to \C$ of functions and a 
$\sigma$-algebra $\B$, we define the \emph{energy} $\Energy_{f}(\B)$ to be the quantity
\begin{equation}\label{energy}
\Energy_{f}(\B) :=\sum_{j=1}^m \| \E(f_j|\B) \|_{L^2}^2.
\end{equation}
\end{definition}

In practice $m$ will either be 1 or 3.  Observe that we have the trivial bounds
\begin{equation}\label{energy-trivial}
0 \leq \Energy_{f}(\B) \leq \sum_{j=1}^m \| f_j \|_{L^2}^2.
\end{equation}
Also from Pythagoras's theorem and the orthogonality considerations discussed above we see that if $\B'$ is finer than $\B$, then
\begin{equation}\label{pythagoras}
 \sum_{j=1}^m \| \E(f_j|\B') - \E(f_j|\B) \|_{L^2}^2 = \Energy_{f}(\B') - \Energy_{f}(\B).
\end{equation}
In particular, the energy of $\B'$ is larger than or equal to $\B$.

We now describe, in abstract terms, the idea of the energy increment strategy.  Suppose one is trying to prove 
a statement $P(M)$ involving some large parameter $M > 0$ which one hopes to keep under control; for instance, one may
be trying to bound some fixed expression $E$ from above by $M$ or from below by $1/M$.  To begin with, this statement
does not depend on any $\sigma$-algebras.  But now we introduce a $\sigma$-algebra $\B$, which we initialize to be
the trivial algebra $\B = \{\emptyset, \Z_N\}$, and try to prove $P(M)$ using an argument which is in some sense ``relative to $\B$'' (in particular, the bounds $M$ may depend on some measure of how ``complex'' $\B$ is).  Either this argument works, or it encounters some obstruction.  The idea is then to show that
the obstruction forces the existence of a new $\sigma$-algebra $\B'$ which is finer than $\B$ (and typically more complex than $\B$) and has slightly more energy.  One then replaces $\B$ by $\B'$ and then repeats the above strategy, hoping to use
the trivial bound \eqref{energy-trivial} to show that the argument must eventually work relative to \emph{some} $\sigma$-algebra.

The difficulty with this strategy is that the energy increment obtained by this method typically depends on 
the complexity of the $\sigma$-algebra $\B$, which tends to grow rather quickly.  As such it is possible for this
method to get bogged down at some intermediate energy in the range \eqref{energy-trivial} and not terminate in any
controlled amount of time.  To get around this, it turns out that instead of just using a pair $\B \subset \B'$ of
$\sigma$-algebras, it is better to use a \emph{triplet} $\B \subset \B' \subset \B''$ of $\sigma$-algebras, with
the energy gap between $\B$ and $\B'$ allowed to be moderately large (bounded by a quantity that does not depend
on the complexity of any of these algebras).  The idea is then to try to prove $P(M)$ relative to the 
pair $(\B,\B')$, but using bounds which depend only on the complexity of $\B$ and not on $\B'$.  If the argument
encounters an obstruction, then one can replace $\B'$ by a more complex $\B''$, with an energy increment again
depending only on the complexity of $\B$; thus this energy increment will not go to zero as $\B'$ becomes more
complex.  There is now a second obstruction when the energy gap $\B'$ and $\B$ becomes too large, but then one
replaces $\B$ by $\B'$; this can only occur a finite number of times because we do not allow the bounds for this energy gap to depend on the complexity and thus the energy increment here is bounded from below by a fixed constant.

To make this argument more precise we encapsulate it in the following abstract lemma (which has a certain resemblance
to Zorn's lemma, and can be in fact thought of as a ``quantitative'' version of that lemma; it also resembles the proof of the Szemer\'edi regularity lemma).  We are indebted to Ben Green for suggesting the use this type of energy incrementation argument, which is for instance used in our joint paper \cite{gt-primes} to establish arbitrarily long arithmetic progressions in the primes.

\begin{lemma}[Abstract energy incrementation argument]\label{abstract}  Suppose there is a property $P(M)$ which can depend on some parameter $M > 0$.  Let $d \geq 0$, and let $f = (f_1,\ldots, f_m)$ be a collection
of $m$ bounded functions.

Suppose also that we have an $\tau > 0$ for which the following dichotomy holds: 
for any $X, X' > 0$, and given any $\sigma$-algebra $\B$ which is compact of order $d$ with complexity at most $X$, and any $\sigma$-algebra $\B'$ which is finer than $\B$ and also compact of order $d$ with complexity at most $X'$, then if the energy gap condition
\begin{equation}\label{egap}
 \Energy_f(\B') - \Energy_f(\B) \leq \tau^2,
\end{equation}
holds, then either $P(M)$ is true for some $M = O_{d,\tau,X,X'}(1)$, or we can find a $\sigma$-algebra $\B''$ finer than $\B'$ which is compact of order $d$ with complexity at most $O_{d,\tau,X,X'}(1)$ such that 
we have the energy increment property
\begin{equation}\label{ei}
 \Energy_f(\B'') - \Energy_f(\B') \geq c(d, \tau,X) > 0
\end{equation}
for some positive quantity $c(d,\tau,X) > 0$ which does \emph{not} depend on $X'$.  

Then $P(M)$ is true for some $M = O_{m,d,\tau}(1)$.
\end{lemma}

\begin{remark}  The point of this lemma is that it reduces the task of proving some property $P(M)$ to the
easier task of proving a \emph{dichotomy}; either $P$ can be proven, or we can increment the energy of a certain $\sigma$-algebra while keeping the complexity under control.  It is crucial that $\tau$ does not depend on the complexities
$X, X'$, and that the energy increment $c(d,\tau,X)$ depends only on the lower complexity $X$ and not the higher complexity $X'$, otherwise this lemma fails.
Note that no quantitative knowledge on the growth of this complexity
on of the energy increment bound $c(d,\tau,X)$ is necessary, although of course the explicit form of the final bound $O_{m,d,\tau}(1)$ on $M$ will depend quite heavily on those growth rates.  
This argument proceeds by a double iteration and thus typically produces bounds which are of Ackermann type or worse,
but in principle they are computable.
\end{remark}

\begin{proof}  The proof proceeds by running the following double iteration algorithm, constructing a 
pair of $\sigma$-algebras
$\B$ and $\B'$, both compact of order $d$ and with $\B'$ finer than $\B$, as follows.

\begin{itemize}

\item[Step 1] Initialize $\B$ to the trivial algebra $\B := \{\emptyset, \Z_N \}$.

\item[Step 2] Initialize $\B'$ to equal $\B$ (thus trivially verifying the energy gap condition \eqref{egap}).
Let $X$ denote the complexity of $\B$.

\item[Step 3] Let $X'$ denote the complexity of $\B'$.  If $P(M)$ is true for some $M = O_{d,\tau,X,X'}(1)$ then we
halt the algorithm.  Otherwise, we must by hypothesis be able to locate a $\sigma$-algebra $\B''$ which is compact of order $d$ with complexity at most $O_{d,\tau,X,X'}(1)$ with the energy increment property \eqref{ei}, and we continue on
to Step 4.

\item[Step 4] If $\Energy_f(\B'') - \Energy_f(\B) \leq \tau^2$, then we replace $\B'$ by $\B''$ (thus preserving
\eqref{egap}) and return to Step 3.  Otherwise, we replace $\B$ by $\B''$ and return to Step 2.

\end{itemize}

Observe that for each fixed $\B$ of complexity $X$, the algorithm can only iterate for at most $O_{d,\tau,X}(1)$
times before changing $\B$.  This is because every time $\B'$ is changed, the energy $\Energy_f(\B')$ increases by
at least $c(d,\tau,X)$, but if the energy ever exceeds $\Energy_f(\B) + \tau^2$ then we must change $\B$.  Note
that it is crucial here that the energy increment $c(d,\tau,X)$ not depend on the complexity $X'$ of $\B'$, which
may be growing quite rapidly during this iteration process.  In particular, if $\B$ finally does change, its complexity
will increase from $X$ to at most $O_{d,\tau,X}(1)$.  Next, observe that $\B$ can only be changed at most
$O_{m,\tau}(1)$ times, because each time we change $\B$, the energy $\Energy_f(\B)$ increases by at least $\tau^2$,
but the energy is always non-negative and is bounded by $m$.  Combining these two observations we see that
the entire algorithm must halt in $O_{m,d,\tau}(1)$ steps and all $\sigma$-algebras constructed by the algorithm
have complexity at most $O_{m,d,\tau}(1)$.  The claim follows.
\end{proof}

\section{Proof of the structure theorem}\label{structure-sec}

We now prove the structure theorem, Theorem \ref{structure}.  Naively, the idea would be to take $\B$ to
be the $\sigma$-algebra formed by \emph{all} the $UAP^{d-1}$ functions, and then take $f_U = f - \E(f|\B)$
and $f_{U^\perp} = \E(f|\B)$ to conclude the result (the uniformity of $f_U$ arising from Proposition \ref{nonrandom}, and Proposition \ref{compact-sigma} being used to locate $f_{UAP}$); this would be the exact analogue of how one would proceed in the genuinely ergodic setting when the underlying space is infinite and one does not care about quantitative control on the complexity of the $\sigma$-algebra.  Unfortunately
this approach does not work because there are far too many $UAP^{d-1}$ functions available, and the complexity of $\B$
would explode with $N$ (indeed, it is likely that $\B$ would simply be the total $\sigma$-algebra consisting of arbitrary subsets of $\Z_N$).  Thus, in the quantitative setting, one must be substantially more ``choosy'' about 
which $UAP^{d-1}$ functions to admit into the algebra $\B$ - they should only be the ones which have a good reason 
for being there, such as having a non-trivial correlation with the function $f$.  It turns out that the best
framework for doing this is given by the abstract energy incrementation argument given in the previous section, exploiting the fact that each function that one adds to the $\sigma$ algebra increases the energy of that algebra, especially if there is a correlation with $f$.  As such, the proof this theorem does not actually require one to know what the function $c_0(d,k,\delta,M)$ in Theorem \ref{recurrence}
actually is (although this function will of course influence the final bound on $M$), and so we can prove this
theorem before proving the (somewhat more difficult) recurrence theorem, Theorem \ref{recurrence}, in the
next section.

In view of the energy incrementation argument, it suffices to prove the following dichotomy:

\begin{lemma}[Structure theorem dichotomy]\label{structure-dich}  Let $k \geq 3$, and let $f$ be a non-negative bounded
function obeying \eqref{f-density} for some $\delta > 0$.  Let $\B \subset \B'$ be $\sigma$-algebras
which are compact of order $k-2$ with complexity at most $X$, $X'$ respectively, and and obey the
energy gap condition \eqref{egap} with $\tau := \frac{\delta^2}{5000 k}$.
Then at least one of the following must be true:

\begin{itemize}

\item (Success) We can find a positive number $M = O_{k,\delta,X}(1)$
a bounded function $f_U$, and non-negative bounded functions $f_{U^\perp}$, $f_{UAP}$ such that we have
the splitting $f = f_U + f_{U^\perp}$ and the estimates \eqref{central-0}, \eqref{central-1}, \eqref{central-2} 
with $d := k-2$, as well as the Gowers uniformity estimate \eqref{u-est}.

\item (Energy increment) We can find a $\sigma$-algebra $\B''$ finer than $\B'$ which is compact of order $d$
and complexity $O_{k,\delta,X,X'}(1)$ such that
\begin{equation}\label{energy-inc}
 \Energy_f(\B'') - \Energy_f(\B') \geq c(k,\delta,X) > 0
\end{equation}
for some $c(k,\delta,X) > 0$ independent of $X'$.

\end{itemize}

\end{lemma}

Indeed, Theorem \ref{structure} follows immediately by applying Lemma \ref{structure-dich} to
Lemma \ref{abstract} (using $m=1$ and using the bounded function $f$, and $\tau := \frac{\delta^2}{5000 k}$).

\begin{proof}[Proof of Lemma \ref{structure-dich}]  Fix $\B$, $\B'$.
Since $\E(f|\B)$ is non-negative and bounded, and $\B$ is compact of order $k-2$ with complexity $O(X)$, 
We may apply Proposition \ref{compact-sigma} to find a non-negative bounded function $f_{UAP}$ such that
\begin{equation}\label{bap}
 \| \E(f|\B) - f_{UAP} \|_{L^2} \leq \frac{\delta^2}{5000 k}
\end{equation}
and
$$ \| f_{UAP} \|_{UAP^{k-2}} < M$$
for some $M = O_{k, \delta, X}(1)$, which we now fix.  From \eqref{bap} and Cauchy-Schwarz we observe that
$$
|\E(f|\Z_N) - \E(f_{UAP}|\Z_N)| = |\E( \E(f|\B) - f_{UAP}|\Z_N)| \leq \frac{\delta^2}{5000 k}$$
and in particular (by \eqref{f-density})
$$ \E(f_{UAP}|\Z_N) \geq \delta/2.$$
Now split $f = f_U + f_{U^\perp}$,
where $f_{U^\perp} := \E(f|\B')$ and $f_U := f - \E(f | \B')$.  We have already proven the estimates
\eqref{central-1}, \eqref{central-2}, while \eqref{central-0} follows 
from \eqref{egap} (recall $\tau = \frac{\delta^2}{5000 k}$), \eqref{pythagoras}, \eqref{bap}, and
the triangle inequality.  If the estimate \eqref{u-est} held then we would now be done (in the ``Success''
half of the dichotomy), so suppose instead that
$$ \| f_U \|_{U^{k-1}} > 2^{-k} c_0(k-2,k,\delta,M).$$
By Lemma \ref{nonrandom}, we can thus find a function $G \in UAP^{d-2}$ with
$ \| G \|_{UAP^{d-2}} \leq 1$
such that
\begin{equation}\label{fug}
 |\langle f_U, G \rangle| > c(k, \delta, M) > 0
\end{equation}
for some positive quantity $c(k,\delta,M)$.  
Now write $\B'' := \B' \vee \B_\eps(G)$, where $\eps = \eps(k, \delta, M) > 0$ is to be chosen later.
We thus split $f_U = (f - \E(f|\B'')) + (\E(f|\B'') - \E(f|\B'))$ and 
$G = (G - \E(G|\B'')) - (\E(G|\B'')-\E(G|\B')) + \E(G|\B')$.  The first terms in both expansions are
orthogonal to $\B''$ (and thus to $\B'$), while the second terms are measurable in $\B''$ and orthogonal to $\B'$,
while the third term of $G$ is measurable in $\B'$.  Thus
$$ \langle f_U, G \rangle = \langle f - \E(f|\B''), G - \E(G|\B'') \rangle
+ \langle \E(f|\B'') - \E(f|\B'), \E(G|\B'') - \E(G|\B') \rangle.$$
From \eqref{g-approx} and the boundedness of $f$ we have
$$ |\langle f - \E(f|\B''), G - \E(G|\B'') \rangle| \leq O(\eps).$$
Thus if we choose $\eps$ sufficiently small depending on $k$, $\delta$, $M$, we see from
\eqref{fug} that
$$ |\langle \E(f|\B') - \E(f|\B), \E(G|\B') - \E(G|\B) \rangle| \geq \frac{1}{2} c(k,\delta,M).$$
Since $G$ is bounded, we thus see from Cauchy-Schwarz that
$$ \| \E(f|\B') - \E(f|\B) \|_{L^2} \geq \frac{1}{4} c(k,\delta,M) > 0.$$
But then \eqref{energy-inc} follows from \eqref{pythagoras}.
Finally, the complexity bound on $\B''$ follows from Definition \ref{complexity-def}, the complexity
bound on $\B$, and the choice of $\eps$ and $M$.  We are thus in the ``Energy increment'' half of the dichotomy,
and the lemma follows.
\end{proof}

The proof of the structure theorem is now complete.

\begin{remark} It may be possible to prove this structure theorem more directly, without explicitly invoking
$\sigma$-algebras, for instance by setting up a extremization problem such as that of minimizing the $U^{k-1}$
norm of $f_U$ subject to the constraints \eqref{central-0}, \eqref{central-1}, \eqref{central-2}, the splitting
$f = f_U + f_{U^\perp}$, and the bounded non-negativity of $f_{UAP}$ and $f_{U^\perp}$.  We were unable however
to achieve this in a clean way, especially when it came to maintaining the boundedness and non-negativity conditions,
whereas the conditional expectation method achieves this more painlessly.
\end{remark}

\section{Compactness on atoms, and an application of van der Waerden's theorem}

To prove Szemer\'edi's theorem, the only thing that now remains is to prove the recurrence theorem for almost
periodic functions, Theorem \ref{recurrence}.  In this section we present a key Proposition, which illustrates the
applicability of van der Waerden's theorem (Theorem \ref{vwt}) to the problem of obtaining recurrence for a function $f$ whose shifts $T^n f$ enjoy a representation such as \eqref{representation}.  The key idea is that
the functions on the right-hand side of \eqref{representation} live in a sufficiently ``compact'' space of functions
that they can be ``finitely coloured'', at which point van der Waerden's theorem can be used to establish
recurrence\footnote{This argument was inspired, not by the original ergodic theory arguments of Furstenberg, but of the later colouring-based arguments, for instance in \cite{bergelson-leibman}.  It may be possible to adapt the older arguments in, say, \cite{furstenberg} instead here, which have the advantage of using the same length $k$ for the progression throughout the argument, instead of replacing $k$ by a considerably larger $k_*$ as is done here.  This might ultimately lead to somewhat better final bounds, although it still seems that one would still get Ackermann-type dependence or worse on the constants.}.  As we show at the end of this section, we can quickly use this Proposition to deduce
the $d=1$ case (as well as the rather trivial $d=0$ case) of Theorem \ref{recurrence} as a corollary.  

\begin{proposition}[Recurrence for conditionally UAP functions]\label{cap}
Let $\B$ be a $\sigma$-algebra, let $M > 0$, let $H$ be a finite non-empty set, and let for
each $n \in \Z_N$ and $h \in H$ let $c_{n,h}$ be a bounded $\B$-measurable function and let $g_h$ be a bounded
function.  Let $h$ be a random variable taking values in $H$, and define the functions $F_n$ for all $n \in \Z_N$ by
the formula
\begin{equation}\label{nrep}
 F_n := M \E( c_{n,h} g_h )
\end{equation}
(compare with \eqref{representation}).  Let $f_{U^\perp}$ be a bounded non-negative function, and for any $\delta > 0$, $n \in \Z_N$, and $k, k_* \in \Z_+$ 
let $E_\lambda(k, \delta, k_*,\B) \in \B$ be the set
\begin{equation}\label{el-def}
\begin{split}
 E_n(k, \delta, k_*,\B) := \{ & x \in \Z_N: \E(T^n f_{U^\perp}|\B)(x) \geq \frac{\delta}{2} \hbox{ and }\\
&\E(|T^n f_{U^\perp} - F_{\lambda m}| | \B)(x) \leq \frac{\delta}{8k} \}.
\end{split}
\end{equation}
Then for every $\delta > 0$ and $k \in \Z_+$ there exists $k_* = k_*(k,\delta,M)$ such that
\begin{equation}\label{eeks} 
\begin{split}
\E( \prod_{j=0}^{k-1} T^{\mu jr} f_{U^\perp}(x)& | x \in \Z_N, 1 \leq r \leq N_0 )\\
&\geq c(k,\delta,M) \E(\P( \bigcap_{m=1}^{k_*} E_{\mu \lambda m}(k, \delta, k_*, \B ) | \Z_N ) | 1 \leq \lambda \leq N_0/k_*)
\end{split}
\end{equation}
for all $\mu \in \Z_N$ and $N_0 \geq k_*$, and some $c(k,\delta) > 0$.
\end{proposition}

\begin{remark}  The point is that this theorem reduces the task of establishing lower bounds for recurrence
expressions involving $f_{U^\perp}$, to that of establishing lower bounds for the recurrence behaviour
of $\B$-measurable sets $E_{\mu \lambda m}(k,\delta,k_*,\B)$.  This is advantageous if
$\B$ is ``simpler'' than the original function $f_{U^\perp}$; in practice, $f_{U^\perp}$ will be approximately an almost periodic function
of order $d$, and $\B$ will be a compact $\sigma$ algebra of order $d-1$, and thus functions measurable in $\B$ can
be approximated by almost periodic functions of one lower order than $d$.  This is the key to the proof of Theorem \ref{recurrence} we give in the next section, which proceeds by induction on $d$.  On the other hand, the bounds on
$k_*$ given by our proof involve van der Waerden numbers, which will cause Ackermann type growth rates or worse
in our final bound.  
\end{remark}

\begin{proof}[Proof of Proposition \ref{cap}]  To prove \eqref{eeks} it suffices to prove the ``localized'' version
\begin{equation}\label{eeks-2}
 \E( \prod_{j=0}^{k-1} T^{\mu\lambda (a+js)} f_{U^\perp}(x) | x \in \Z_N; 1 \leq a, s \leq k_*/k )
\geq c(\delta,k,k_*) \P( \bigcap_{m=1}^{k_*} E_{\mu\lambda m}(k, \delta, k_*,\B ) | \Z_N ) 
\end{equation}
for each $1 \leq \lambda \leq N_0/k_*$.  Indeed, if \eqref{eeks-2} held then upon averaging in $\lambda$ we obtain
\begin{align*}
 \E( \E(\prod_{j=0}^{k-1} T^{\mu \lambda (a+js)} f_{U^\perp}|\Z_N) | &1 \leq \lambda \leq N_0/k_*; 1 \leq a, s \leq k_*/k )\\
&\geq c(\delta,k,k_*) \E( \P( \bigcap_{m=1}^{k_*} E_{\mu\lambda m}(k, \delta, k_*,\B ) | \Z_N ) | 1 \leq \lambda \leq N_0/k_* ).
\end{align*}
The $T^{\mu\lambda a}$ can be factored out of the product and makes no difference to the expectation, thus
it can be discarded.  The $a$ averaging then becomes redundant, and we obtain
\begin{align*}
 \E( \E(\prod_{j=0}^{k-1} T^{\mu \lambda js} f_{U^\perp}|\Z_N) | &1 \leq \lambda \leq N_0/k_*; 1 \leq s \leq k_*/k )\\
&\geq c(\delta,k,k_*) \E( \P( \bigcap_{m=1}^{k_*} E_{\mu\lambda m}(k, \delta, k_*,\B ) | \Z_N ) | 1 \leq \lambda \leq N_0/k_* ).
\end{align*}
The claim \eqref{eeks} then follows by observing that every $1 \leq r \leq N_0$ has at most $O_{k_*,k}(1)$ representations 
of the form $r = \lambda (a+js)$ with  $1 \leq \lambda \leq N_0/k_*$ and $1 \leq a,s \leq k_*/k$.  (The dependence
of $k_*$ is ultimately irrelevant since $k_*$ itself will ultimately depend on $\delta, k, M$).  

It remains to prove \eqref{eeks-2}.
Fix $\mu$, $\lambda$.  By absorbing $\mu$ into $\lambda$ we may take $\mu=1$ (we will not use the upper or lower bound
on $\lambda$). Since $\bigcap_{m=1}^{k_*} E_{\lambda m}(k,\delta,k_*,\B)$ is measurable in $\B$, it is the union of atoms $A \in \B$.  It will suffice to prove the pointwise estimate
$$ \E( \E(\prod_{j=0}^{k-1} T^{\lambda (a+js)} f_{U^\perp}|A) | 1 \leq a,s \leq k_*/k ) \geq c(\delta,k,k_*)$$
for each such atom, as the claim then follows by multiplying this formula by $\P(A|\Z_N)$ and summing over all
atoms in $E_\lambda(k,\delta,k_*,\B)$.

Now fix the atom $A$.  Since the number of pairs $(a,s)$ is $O_{k,k_*}(1)$, it suffices to locate a \emph{single}
pair $(a,s)$ with $1 \leq a,s \leq k_*/k$ such that
\begin{equation}\label{shifty}
 \E(\prod_{j=0}^{k-1} T^{\lambda (a+js)} f_{U^\perp}|A) \geq c(\delta,k)
\end{equation}
for some $c(\delta,k) > 0$.

We now pass from the shifts $T^{\lambda (a+js)} f_{U^\perp}$ to the functions $F_{\lambda (a+js)}$.  
We claim that to prove \eqref{shifty} it would suffice to prove that
\begin{equation}\label{shifty-2}
 \| F_{\lambda (a+js)} - F_{\lambda a} \|_{L^2(A)} \leq \frac{\delta}{8k} \hbox{ for all } 0 \leq j \leq k-1,
\end{equation}
where $L^2(A)$ is the Hilbert space given by the norm $\|F\|_{L^2(A)} := \E( |F|^2 | A)^{1/2}$.
To see this claim, observe from Cauchy-Schwarz that \eqref{shifty-2} implies
$$ \E( |F_{\lambda (a+js)} - F_{\lambda a}| \bigr| A ) \leq \frac{\delta}{8k} \hbox{ for all } 0 \leq j \leq k-1.$$
But by \eqref{el-def} and the choice of $A$, we also have
$$ \E( |F_{\lambda (a+js)} - T^{\lambda (a+js)} f_{U^\perp}| \bigr| A ) \leq \frac{\delta}{8k} \hbox{ for all } 0 \leq j \leq k-1.$$
From the triangle inequality we thus have
$$ \E( |T^{\lambda(a+js)} f - T^{\lambda a} f_{U^\perp}| \bigr| A ) \leq \frac{3\delta}{8k} \hbox{ for all } 0 \leq j \leq k-1.$$
This in particular implies
$$ \E( {\bf 1}_{|T^{\lambda(a+js)} f_{U^\perp} - T^{\lambda a} f_{U^\perp}| > \frac{4}{5} T^{\lambda a} f} T^{\lambda a} f | A)
\leq \frac{15\delta}{32 k} \hbox{ for all } 0 \leq j \leq k-1.$$
On the other hand, by \eqref{el-def} again we have
$$ \E( T^{\lambda a} f_{U^\perp} | A ) \geq \delta/2.$$
Thus
$$ \E( {\bf 1}_{|T^{\lambda(a+js)} f_{U^\perp} - T^{\lambda a} f_{U^\perp}| \leq \frac{4}{5} T^{\lambda a} f
\hbox{ for all } 0 \leq j \leq k-1} T^{\lambda a} f_{U^\perp} )
\geq \frac{\delta}{32 k}.$$
By H\"older's inequality and the non-negativity of $f_{U^\perp}$ this implies that
$$ \E( {\bf 1}_{|T^{\lambda(a+js)} f_{U^\perp} - T^{\lambda a} f_{U^\perp}| \leq \frac{4}{5} T^{\lambda a} f_{U^\perp}
\hbox{ for all } 0 \leq j \leq k-1} (T^{\lambda a} f_{U^\perp})^k )
\geq (\frac{\delta}{32 k})^k.$$
The claim \eqref{shifty} then follows from the elementary pointwise inequality
$$ \prod_{j=0}^{k-1} T^{\lambda (a+js)} f_{U^\perp} \geq \frac{1}{5}^k
{\bf 1}_{|T^{\lambda(a+js)} f_{U^\perp} - T^{\lambda a} f_{U^\perp}| \leq \frac{4}{5} T^{\lambda a} f_{U^\perp}
\hbox{ for all } 0 \leq j \leq k-1} (T^{\lambda a} f_{U^\perp})^k.$$

It remains to find a pair $(a,s)$ obeying \eqref{shifty-2}.  Using \eqref{nrep} it suffices to find an $(a,s)$
such that
\begin{equation}\label{shifty-3}
\| \E( c_{\lambda(a+js),h} g_h ) - \E( c_{\lambda a, h} g_h ) \|_{L^2(A)} \leq \delta/8Mk
\hbox{ for all } 0 \leq j \leq k-1.
\end{equation}
Note that as the $c_{n,h}$ are measurable with respect to $\B$, they are constant on $A$, and so
without loss of generality we can treat them just as bounded complex numbers (this is the whole point of
working on individual atoms in the first place).  The $g_h$ are not constant, but we can think of
them as bounded functions on $A$.

To proceed further we need the following compactness property of averages of the form $\E( c_h g_h )$
in $L^2(A)$.

\begin{lemma}[Total boundedness property]\label{tbp}  
There exists integers $1 \leq m_1, \ldots, m_L \leq k_*$ for some $L \leq C(k,M,\delta)$
such that
$$ \inf_{1 \leq l \leq L} \| \E( c_{\lambda m, h} g_h ) - \E( c_{\lambda m_l, h} g_h ) \|_{L^2(A)} 
\leq \frac{\delta}{16Mk} \hbox{ for all } 1 \leq m \leq k_*.$$
\end{lemma}

\begin{remark} The key point here is that the bound on $L$ does not depend on the size of $H$, $A$, or $N$.
This is a quantitative analogue of the basic result (used in the ergodic theory proofs, see e.g. \cite{furstenberg}) that
a Volterra integral operator from one finite measure space to another is necessarily a compact operator in $L^2$, and thus
the range of any bounded set can be covered by a finite number of $\delta$-balls in $L^2$.
\end{remark}

\begin{proof}  Let us write $f_m := \E( c_{m,h} g_h )|_A$.
We construct an orthonormal system of functions $v_1, v_2, \ldots, v_J$ in $L^2(A)$ by performing
the following algorithm, which can be viewed as a rudimentary version of the energy increment algorithm discussed
in previous sections (with the role of $\sigma$-algebras replaced by the simpler notion of finite-dimensional subspaces
of a Hilbert space).

\begin{itemize}

\item[Step 0] Initialize $J=0$.

\item[Step 1] Let $V \subset L^2(A)$ be the subspace spanned by the $v_1, \ldots, v_J$ (so initially this will
be the trivial space $\{0\}$).  

\item[Step 2] If there exists a $1 \leq m \leq k_*$ such that $\dist_{L^2(A)}(f_m, V) \geq \delta/64Mk$, then by 
Hilbert space geometry
we can find a unit vector $v_{J+1}$ orthogonal to $V$ (and thus to all the $v_1, \ldots, v_J$) such that
$|\langle f_m, v_{J+1} \rangle_{L^2(A)}| \geq \delta/64Mk$.  In such a case, we choose\footnote{Note that since $m$ ranges over a finite set, the axiom of choice is not needed here, since $\Z_N$ is clearly well-ordered.  Because we are always in a finite (or at least finite dimensional) setting, similar considerations apply to other parts of the argument in which an arbitrary choice has to be made.} such a $v_{J+1}$, increment $J$, and
return to Step 1.  Otherwise, we terminate the algorithm.

\end{itemize}

We claim that this algorithm terminates in $O_{k,M,\delta}(1)$ steps.  Indeed, for each $v_j$ generated by this algorithm,
we see from construction that there exists an $m = m(j) \in \Z_N$ such that
$$ |\E( c_{\lambda m,h} \langle g_h, v_j \rangle_{L^2(A)} )| = |\langle f_m, v_j \rangle_{L^2(A)}| \geq \delta/64Mk.$$
Here we have crucially taken advantage of the fact that $c_{\lambda m, h}$ is constant on $A$.
Since $c_{\lambda m,h}$ is bounded, we thus see from the Cauchy-Schwarz inequality that
$$ \E( |\langle g_h, v_j \rangle_{L^2(A)}|^2 ) \geq (\frac{\delta}{64 M k})^2.$$
Summing this in $j$, we obtain
$$ \E( \sum_{j=1}^J |\langle g_h, v_j \rangle_{L^2(A)}|^2 ) \geq (\frac{\delta}{64 M k})^2 J.$$
But from the boundedness of the $g_h$, the orthonormality of the $v_j$, and Bessel's inequality, the left-hand side
is at most 1.  Thus $J \leq (\frac{64 Mk}{\delta})^2 = O_{k,M,\delta}(1)$ as claimed.

Now observe from the construction of the algorithm that all the functions $f_m$ will lie within
$\delta/64Mk$ (in the $L^2(A)$ metric) of the $J$-dimensional space $V$.  In particular, we see from the triangle inequality, the crude bound $\|f_m\|_{L^2(A)} \leq 1$ arising from our bounds on $c_{n,h}$ and $g_h$, and finite-dimensional  geometry that there can be at most $O_{k,\delta,J}(1) = O_{k,M,\delta}(1)$ functions $f_{m_1}, \ldots, f_{m_L}$ which are all separated from each other by at least $\delta/16Mk$ in the $L^2(A)$ metric.  The claim now follows 
by the usual greedy algorithm.
\end{proof}

Using this lemma, we can introduce a colouring function $\c: \{1,\ldots,k_*\} \to \{1,\ldots,L\}$
by
$$ \c(m) := \inf \{ 1 \leq l \leq L: \| \E( c_{\lambda m, h} g_h ) - \E( c_{\lambda m_l, h} g_h ) \|_{L^2(A)} 
\leq \delta/16Mk \}.$$
By van der Waerden's theorem, if $k_* = k_*(k,L) = k_*(k,\delta,M)$ is chosen sufficiently large, then
we can find $1 \leq a, s \leq k_*/k$ such that the progression $a, a+s, \ldots, a+(k-1)s$ is monochromatic.  The claim
\eqref{shifty-3} now follows from the triangle inequality.  This concludes the proof of Proposition \ref{cap}.
\end{proof}

As a quick corollary of this Proposition we can now prove the $d=1$ case, at least, of Theorem \ref{recurrence}.

\begin{proof}[Proof of Theorem \ref{recurrence} when $d=1$]
Let $f_{U^\perp}$, $f_{UAP}$, $k$, $M$, $\delta$, $\epsilon$ be as in the Theorem.  From \eqref{central-2} and Definition \ref{ap-def} we can find a finite non-empty set $H$, a collection of bounded constants $(c_{n,h})_{n \in \Z_N; h \in H}$, and bounded functions $(g_h)_{h \in H}$, and a random variable $h$ taking values in $H$ such that we have the
representation \eqref{nrep}, where $F_n := T^n f_{UAP}$.  We thus apply Proposition \ref{cap} with 
$N_0 := N_1$ and $\B$ set equal to the trivial $\sigma$-algebra $\B = \{\emptyset, \Z_N\}$, since the $c_{n,h}$ are all almost periodic of order 0 and hence constant.  But by \eqref{el-def}, we see that $
E_{\lambda m}(k,\delta,k_*,\B)$ is either
the empty set or all of $\Z_N$, with the latter occuring if
$$ \E(T^{\lambda m} f_{U^\perp}|\Z_N) \geq \delta/2 \hbox{ and } \E(|T^{\lambda m} f - F_{\lambda m}| | \Z_N) \leq \frac{\delta}{8k}.$$
But the latter condition is automatic from \eqref{central-1}, while the latter follows from \eqref{central-0},
Cauchy-Schwarz, and the choice of $\epsilon$; note that the shift $T^{\lambda m}$ has no effect
on the expectation $\E(|\Z_N)$.  Thus $E_{\lambda m}(k,\delta,k_*,\B) = \Z_N$ for all $\lambda$,
and the claim \eqref{cdm} follows from \eqref{eeks}.
\end{proof}

\begin{remark}
As we shall see, the $d > 1$ case is somewhat more complicated, the problem beign that one has to somehow ``quotient out'' the effect of a very large number of almost periodic functions of order $d-1$ before the property of being almost periodic of order $d$ emerges as a usable property.  This appears to unfortunately be rather
necessary, even when $d=2$, at least with the arguments currently available; the author would consider this issue of the least well understood components of the theory.  Consider for instance a function $f_{AP}$ of the form $f_{AP}(x) = \psi( x^2 / N )$, where $\psi: \R/Z \to [0,1]$ is a smooth bounded non-negative function which 
is periodic with period 1, which equals 1 on the interval $[-\delta,\delta]$, and vanishes outside $[-2\delta,2\delta]$.
This function can be shown to be almost periodic of order 2 with an $UAP^2$ norm of $O_\delta(1)$.  Thus Theorem
\ref{recurrence} should allow us to locate a large number of arithmetic progressions of length $k$ in the support
of $\psi$, for reasonably large values of $k$ (e.g. $k=5$).  To actually establish even this special case, however, seems rather difficult, the simplest proof probably being the ergodic theory proof that lifts this problem up to establishing
recurrence for the skew shift on the two-dimensional torus.  Similarly for more complicated examples such as \eqref{fab} (now the ergodic system is a two-step nilsystem, formed by quotienting the unipotent upper triangular $3 \times 3$
matrices by the subgroup of matrices with integer coefficients).  In \cite{gowers} this precise problem was encountered, and solved by using very directly the number-theoretic structure of $x^2/N$ (and similarly polynomial objects), in particular a quantitative version of Weyl's theorem on the uniform distribution of polynomials.  The problem of
having to deal with generalized polynomials instead of polynomials was avoided by working on relatively short
arithmetic progressions, in which one could approximate the former by the latter.
\end{remark}

\section{Recurrence for almost periodic functions}\label{central-sec}

We now conclude the proof of Theorem \ref{recurrence}, and thus of Theorem \ref{szt}.  We have already
handled the $d=1$ case.  The case $d=0$ can either be deduced from the $d=1$ case, or can be worked out directly
by an easy argument which we leave to the reader (the point being that $f_{UAP}$ is now constant and $f_{U^\perp}$
and its shifts will have to be larger than, say, $\delta/2$ with very high probability, say at least $1-1/2k$).
Thus it remains to handle the $d > 1$ cases.  We may assume as an inductive hypothesis that $d$ is fixed
and the claim has already been proven for $d-1$.

When $\mu=0$ the claim follows easily from \eqref{central-1} and the boundedness of $f_{U^\perp}$, so we may
take $\mu \neq 0$.  But then we may rescale by $\mu$ and set $\mu=1$.

We would like to apply Proposition \ref{cap} as we did in the $d=1$ case.  The difficulty now is that the 
functions $c_{n,h}$ generated by Definition \ref{ap-def} are no longer constant, but are themselves almost
periodic of one lower order, $d-1$.  The strategy is then to locate a $\sigma$-algebra $\B$ generated by such functions
(and hence compact of order $d-1$) with respect to which the $c_{n,h}$ are close to being measurable (i.e. close
to constant on most atoms).  Proposition
\ref{cap} then allows us to reduce the problem of establishing recurrence for $f_{U^\perp}$ to one of establishing
a property very similar to recurrence for certain subsets of $\B$, which we can handle by combining the induction
hypothesis with Proposition \ref{compact-sigma}.  As with the structure theorem, one would naively want to take $\B$
to be the $\sigma$-algebra generated by \emph{all} the $c_{n,h}$ (and this is indeed what one does in the genuinely ergodic setting), but again we lose control of the complexity this way.
Instead we must be much more selective with which $c_{n,h}$ we admit.  Again, the easiest framework to implement
this idea is given by the abstract energy increment lemma, Lemma \ref{abstract}.
The point is that it may happen that the $c_{n,h}$ are refusing to be close to measurable on $\B$, or that other problems arise such as $\B$ failing to be sufficiently ``shift-invariant'' (this issue arose in the $d=1$ case when one needed
to eliminate the $T^{\lambda m}$ shift, although in that case the resolution to the problem was trivial).  
In that case, however, the simplest solution is to replace $\B$ by a larger
$\sigma$-algebra $\B'$, to which one adds in all the obstructions (or at least a representative sample thereof) which one encountered in closing the argument, thus increasing the energy of $\B$.  

We turn to the details.  It will suffice to establish

\begin{proposition}[Recurrence theorem dichotomy]\label{recurrence-dich}  
Let $d \geq 2$, and suppose that Theorem \ref{recurrence} has already been proven for $d-1$.  Let 
$k \geq 1$ be integers, and let $M, \delta > 0$.  All quantities in
what follows can depend on $d$, $\delta$, $k$, $M$ (including the implicit bounds in $O()$ notation), and we omit future dependence on these parameters.
let $f_{U^\perp}, f_{UAP}$ be non-negative bounded functions obeying the bounds \eqref{central-0}, \eqref{central-1}, \eqref{central-2}.  Write $f := (f_{U^\perp}, |f_{U^\perp} - f_{UAP}|, |f_{U^\perp} - f_{UAP}|^2)$.
Let $\B \subset \B'$ be $\sigma$-algebras
which are compact of order $k-2$ with complexity at most $X$, $X'$ respectively, and such that
\eqref{egap} holds for some small $\tau > 0$ independent of $X$, $X'$
to be chosen later.  Then at least one of the following must be true:

\begin{itemize}

\item (Success) We have
\begin{equation}\label{cdm-extend}
  \E( \prod_{j=0}^{k-1} T^{jr} f_{U^\perp}(x) | x \in \Z_N, 0 \leq r \leq N_1) \geq c(\tau,X)
\end{equation}
for some $c(\tau,X) > 0$ and all $N_1 \geq 0$.

\item (Energy increment) We can find a $\sigma$-algebra $\B''$ finer than $\B'$ which is compact of order $d$
and complexity $O_{\tau,X,X'}(1)$ such that 
\begin{equation}\label{energy-inc-2}
 \Energy_f(\B'') - \Energy_f(\B') \geq c(\tau,X) > 0.
\end{equation}
Note that the increment $c(\tau,X)$ in \eqref{energy-inc-2} is independent of $X'$.

\end{itemize}

\end{proposition}

Indeed, Theorem \ref{recurrence} will follow from this Lemma and Lemma \ref{abstract} (setting $m = 3$
and letting $1/c(\tau,X)$ play the role of $M$).

It remains to prove Proposition \ref{recurrence-dich}.  We will aim towards applying Proposition \ref{cap},
by locating a large subset of $\B'$ where $f_{U^\perp}$ (and several of its shifts) are large on average,
$f_{UAP}$ is close to $f_{U^\perp}$ on average (as are various shifts of these functions), and the $c_{n,h}$
are close to constant, and then using the induction hypothesis to obtain lower bounds on the sets $E_{\lambda\mu}$ obtained this way.  There may be some obstructions to implementing this strategy, but when they arise
we will convert those obstructions to an energy increment, establishing \eqref{energy-inc-2} instead of
\eqref{cdm-extend}.

\begin{proof}[Proof of Lemma \ref{recurrence-dich}]  
By \eqref{central-2} and Definition \ref{ap-def} we can find a 
finite non-empty set $H$, a collection of bounded functions $(c_{n,h})_{n \in \Z_N; h \in H}$ in $UAP^{d-1}$ with
$\|c_{n,h}\|_{UAP^{d-1}} \leq 1$, and bounded functions $(g_h)_{h \in H}$, and a random variable $h$ taking values 
in $H$ such that we have the representation
\begin{equation}\label{fap-n}
 T^n f_{UAP} = M\E( c_{n,h} g_h )
\end{equation}
for all $n \in \Z_N$.  We cannot yet apply Proposition \ref{cap} since the $c_{n,h}$ are not necessarily measurable
with respect to $\B'$; indeed there are too many of the $c_{n,h}$ to safely add all of them to $\B'$, which needs to
have bounded complexity.  Instead, we shall work using much smaller batches of $c_{n,h}$ and then average
at the end.

We will need a large integer $N_0 = N_0(\tau,X)> 1$ to be chosen later\footnote{There are a number of parameters involved here, which
are at several different scales.  In order to have some idea of what parameters should be large and what parameters
should be small, we suggest using the hierarchy
$$ d, k, \frac{1}{\delta}, M, k_* \ll \frac{1}{\tau} \ll X \ll N_0 \ll X' \ll N_1, N, |H|$$
which is a very typical arrangement of the parameters.  The key points
are that the energy gap $\tau$ does not depend on the large parameters $X, N_0, X', N_1, N, |H|$, that the
energy increment in \eqref{energy-inc} does not depend on the very large parameters $X', N_1, N, |H|$, and
the remaining bounds do not depend on the extremely large parameters $N_1, N, |H|$.}
If $N_1 \leq N_0$ then the claim \eqref{cdm-extend}
follows easily from \eqref{central-1} just by considering the $r=0$ component of the left-hand side, so we will
assume $N_1 > N_0$.  We then observe that
\begin{equation}\label{avast}
 \E( \prod_{j=0}^{k-1} T^{jr} f_{U^\perp}(x) | x \in \Z_N; 0 \leq r \leq N_1) \geq c(N_0)
\E( \E( \prod_{j=0}^{k-1} T^{\mu jr} f_{U^\perp}(x) | x \in \Z_N, 1 \leq r \leq N_0) | 1 \leq \mu \leq N_1/N_0 )
\end{equation}
because each $0 \leq r \leq N_0$ has at most $O_{N_0}(1)$ representations of the form $r = \mu r'$ where 
$1 \leq r' \leq N_0$ and $1 \leq \mu \leq N_1/N_0$.

Now fix a single $1 \leq \mu \leq N_1/N_0$, and consider the expression
\begin{equation}\label{mu-form}
 \E( \prod_{j=0}^{k-1} T^{\mu jr} f_{U^\perp}(x) | x \in \Z_N, 1 \leq r \leq N_0).
\end{equation}
Observe that the exponents $\mu j r$ now range in the relatively small set $\mu \cdot \{ 0, \ldots, (k-1)N_0 \}$.
This has localized the ``$n$'' index in \eqref{fap-n} to a reasonably bounded set (one which is independent of $N$), but
the ``$h$'' parameter is still ranging over a potentially unbounded set $H$.  To resolve this we need the following variant
of Lemma \ref{tbp}.

\begin{lemma}[Finite-rank approximation]\label{fra}  Let $\mu \in \Z_N$.  
Then we can find
$h_1, \ldots, h_{N_0^{100}} \in H$ (not necessarily distinct, and depending on $\mu$) such that
\begin{equation}\label{geo}
\| \E(c_{\mu m,h} g_h) - \E(c_{\mu m,h_j} g_{h_j} | 1 \leq j \leq D) \|_{L^2} \leq O(N_0^{-40})
\end{equation}
for all $0 \leq m \leq (k-1)N_0$.
\end{lemma}

\begin{proof}  We use the second moment method.  Set $D := N_0^{100}$, $G_h := c_{\mu m,h} g_h$, $F := \E(G_h)$,
and let $h_1, \ldots, h_D$ be $D$ independent samples of
the random variable $h$.  We will show that
$$ \P( \| F - \E(G_{h_j} | 1 \leq j \leq D) \|_{L^2} > \frac{(kN_0)^{1/2}}{N_0^{50}} ) 
\leq \frac{1}{kN_0},$$
which implies the claim with probability at least $1 - \frac{(k-1)N_0 + 1}{kN_0} > 0$.

By Chebyshev's inequality, it will suffice to show that
$$ \E( \| F - \E(G_{h_j} | 1 \leq j \leq D) \|_{L^2}^2 ) \leq 1/D.$$
The left-hand side can be expanded as
$$ \E( \E( |F(x)|^2 - 2 \Re \E( \overline{F(x)} G_{h_j}(x) | 1 \leq j \leq D)
+ |\E(G_{h_j}(x) | 1 \leq j \leq D)|^2 | x \in \Z_N) )$$
which we expand and rearrange further as
\begin{equation}\label{ff}
 \| F \|_{L^2}^2 - 2 \Re \E( \overline{F(x)} \E(G_{h_j}(x)) | 1 \leq j \leq D; x \in \Z_N )
+ \E( \E(\overline{G_{h_{j'}}(x)} G_{h_j}(x)) | 1 \leq j,j' \leq D; x \in \Z_N ).
\end{equation}
Since $h_j, h_{j'}$ were chosen with the same distribution as $h$, and are independent when $j \neq j'$, we have
$$ \E( G_{h_j}(x) ) = \E( G_h(x) ) = F(x)$$
and
$$ \E( \overline{G_{h_{j'}}(x)} G_{h_j}(x)) ) = \E( \overline{G_{h'}(x)} | h' \in H ) \E( G_h(x) | h \in H )
= |F(x)|^2 \hbox{ when } j \neq j'.$$
We thus can rewrite \eqref{ff} as
$$ \| F \|_{L^2}^2 - 2 \| F \|_{L^2}^2 + \| F \|_{L^2}^2 + \E(
 \delta_{j,j'} (\E(\overline{G_{h_{j'}}(x)} G_{h_j}(x)) - |F(x)|^2) | 1 \leq j,j' \leq D; x \in \Z_N )$$
where $\delta_{j,j'}$ is the Kronecker delta.  When $j=j'$, we have
$$ \E(\overline{G_{h_{j'}}(x)} G_{h_j}(x)) - |F(x)|^2 = \E( |G_h(x)|^2 | h \in H ) - |\E( G_h(x) | h \in H)|^2$$
which is at most $1$ since $G_h$ is bounded (in fact one can sharpen this to $\frac{1}{4}$, but we will not need
this).  The claim follows.
\end{proof}

Let $h_1, \ldots, h_{N_0^{100}}$ be as in the above Lemma.  Then from \eqref{geo} and \eqref{fap-n} we see that
\begin{equation}\label{tmum}
 \| T^{\mu m} f_{UAP} - M\E( c_{\mu m,h_j} g_{h_j} | 1 \leq j \leq N_0^{100} ) \|_{L^2} \leq O(N_0^{-40}) \hbox{
for all } 0 \leq m \leq (k-1)N_0.
\end{equation}
We have now modeled a reasonably large number of shifts of our almost periodic function $f_{UAP}$
in terms of a controlled number of functions $c_{\mu m, h_j}$.  Next, we define a new $\sigma$-algebra $\B''$
finer than $\B'$ (and depending on $\mu$, $h_1, \ldots, h_{N_0^{100}}$) by
$$ \B'' := (\bigvee_{-(k-1)N_0 \leq m \leq (k-1)N_0} T^{\mu m} \B') \vee (\bigvee_{0 \leq m \leq (k-1)N_0; 1 \leq j \leq N_0^{100}} \B_{N_0^{-100}}(c_{\mu m, h_j}))$$
where $\B_\eps(G)$ are the $\sigma$-algebras constructed by Proposition \ref{ap-generate}.
Since the $c_{\mu m,h_j}$ are in $UAP^{d-1}$ with norm at most 1, and $\B'$ was compact of order $d-1$
and complexity at most $X'$, we see from \eqref{shift-invariant} and
Definition \ref{complexity-def} that $\B''$ is also compact of
order $d-1$ and complexity at most $O_{N_0, X'}(1)$.  Also, from \eqref{g-approx} we have
$$ \| c_{\mu m, h_j} - \E( c_{\mu m, h_j} | \B'' ) \|_{L^\infty} \leq O(N_0^{-100}) \hbox{ for all }
0 \leq m \leq (k-1)N_0$$
and hence (since the $g_{h_j}$ are bounded)
$$
 \| M\E( c_{\mu m,h_j} g_{h_j} | 1 \leq j \leq D ) - M\E( \E(c_{\mu m,h_j}|\B'') g_{h_j} | 1 \leq j \leq D )
\|_{L^2} \leq O(N_0^{-100}).
$$
Combining this with \eqref{tmum} we see that
\begin{equation}\label{av-ok}
 \| T^{\mu m} f_{UAP} - F_{\mu m} \|_{L^2} \leq O_{M}(N_0^{-40}) \hbox{ for all } 0 \leq m \leq (k-1)N_0,
\end{equation}
where $F_n$ is defined for $n \in \Z_N$ by the formula
\begin{equation}\label{fn-def}
F_n := M\E( \E(c_{n,h_j}|\B'') g_{h_j} | 1 \leq j \leq D ).
\end{equation}
We may then apply Proposition \ref{cap} (with $\B$ replaced by $\B''$) to estimate the quantity \eqref{mu-form} as
\begin{equation}\label{mu-squat}
 \eqref{mu-form} 
\geq c \E(\P( \bigcap_{m=1}^{k_*} E_{\mu \lambda} | \Z_N ) | 1 \leq \lambda \leq N_0/k_*)
\end{equation}
where $k_* = O(1)$ and $E_{\mu \lambda m} = E_{\mu \lambda m}(k,\delta, k_*,\B'')$ was defined in
that Proposition; recall that we are suppressing all dependence on the quantities $d$, $k$, $\delta$, $M$.

Our attention thus turns to obtaining lower bounds for the size of $E_{\mu \lambda}$.  We first use
\eqref{av-ok} to pass from $F_{\mu m}$ back to $T^{\mu m} f_{UAP}$ (modulo errors that can be made small by making $D$ large).  From \eqref{av-ok} and Cauchy-Schwarz we have
$$ \E( \E( |T^{\mu m} f_{UAP} - F_{\mu m}| | \B'') | \Z_N) = 
\E( |T^{\mu m} f_{UAP} - F_{\mu m}| | \Z_N) \leq O(N_0^{-40}) \hbox{ for all } 0 \leq m \leq (k-1)N_0,$$
so by Markov's inequality
$$ \P( \E( |T^{\mu m} f_{UAP} - F_{\mu m}| | \B'') \geq \frac{\delta}{16k} | \Z_N) \leq
O(N_0^{-40}) \hbox{ for all } 0 \leq m \leq (k-1)N_0.$$
In particular, from \eqref{el-def} and the triangle inequality we see (since $k_* = O(1)$) that
\begin{equation}\label{e-ep}
 \P( \bigcap_{m=1}^{k_*} E_{\mu \lambda m} | \Z_N ) \geq \P( \bigcap_{m=1}^{k_*} E'_{\mu \lambda m} | \Z_N ) - O(N_0^{-30}),
\end{equation}
where
\begin{align*}
 E'_n := \{ &x \in \Z_N: \E(T^n f_{U^\perp}|\B'')(x) \geq \delta/2 \hbox{ and }\\
&\E(|T^n f_{U^\perp} - T^n f_{UAP}| | \B'')(x) \leq \frac{\delta}{16k}\}.
\end{align*}
The next step is to pull the shifts $T^n$ out of the $\B''$ expectations.  To do this we use
the following observation.

\begin{lemma}[Effective shift invariance of $\B''$]  Suppose that $-(k-1)N_0 \leq m \leq (k-1)N_0$ is such that
$$ \| \E( T^{\mu m} f_{U^\perp} | \B'' ) - T^{\mu m} \E( f_{U^\perp} | \B'' ) \|_{L^2} \geq N_0^{-100}$$
or
$$ \| \E( T^{\mu m} |f_{U^\perp}-f_{UAP}|^2 | \B'' ) - T^{\mu m} \E( |f_{U^\perp}-f_{UAP}|^2 | \B'' ) \|_{L^2} \geq N_0^{-100}.$$
Then we are in the energy increment half of the dichotomy of Proposition \ref{recurrence-dich}.
\end{lemma}

\begin{proof} We prove the first claim, as the second is analogous.  Observe that
$\E( T^{\mu m} f_{U^\perp} | \B'' ) = T^{\mu m} \E( f_{U^\perp} | T^{-\mu m} \B'' )$, and so
$$ \| \E( f_{U^\perp} | T^{-\mu m} \B'') - \E( f_{U^\perp} | \B'' ) \|_{L^2} \geq N_0^{-100}.$$
By the triangle inequality, we thus have either
$$ \|  \E( f_{U^\perp} | T^{-\mu m} \B'') - \E( f_{U^\perp} | \B' ) \|_{L^2} \geq \frac{1}{2} N_0^{-100}$$
or
$$ \|  \E( f_{U^\perp} | \B'') - \E( f_{U^\perp} | \B' ) \|_{L^2} \geq \frac{1}{2} N_0^{-100}.$$
But in either case we can use \eqref{pythagoras} (observing that $\B''$ and $T^{-\mu m} \B''$ are both finer than
$\B'$, by construction of $\B''$) to obtain an energy increment \eqref{energy-inc-2}, as desired.  Similarly
for the second claim (which uses the second component of $f = (f_{U^\perp}, |f_{U^\perp} - f_{UAP}|, |f_{U^\perp} - f_{UAP}|^2)$ rather than the first).
\end{proof}

In light of this lemma, we may assume that
$$ \| \E( T^n f_{U^\perp} | \B'' ) - T^n \E( f_{U^\perp} | \B'' ) \|_{L^2}, 
\| \E( T^n |f_{U^\perp}-f_{UAP}|^2 | \B'' ) - T^n \E( |f_{U^\perp}-f_{UAP}|^2 | \B'' ) \|_{L^2} \leq N_0^{-100}$$
for all $n \in \Z_N$.  In particular we have
$$ \P( |\E(T^{\mu m} f_{U^\perp}|\B'')(x) - T^{\mu m} \E(f_{U^\perp}|\B'')(x)| \geq \delta/4 )
\leq O(N_0^{-100})$$
and
$$ \P( |\E(|T^{\mu m} f_{U^\perp} - T^{\mu m} f_{UAP}| |\B'')(x) - T^{\mu m} \E(|f_{U^\perp}-f_{UAP}||\B'')(x)| \geq \frac{\delta}{32k} ) \leq O(N_0^{-50})$$
for all $0 \leq m \leq (k-1) N_0$.  This allows us to estimate (since $k_* = O(1)$)
\begin{equation}\label{ep-epp}
 \P( \bigcap_{m=1}^{k_*} E'_{\mu \lambda m} | \Z_N ) \geq \P( \bigcap_{m=1}^{k_*} E''_{\mu \lambda m} 
| \Z_N ) - O(N_0^{-50}),
\end{equation}
where
$$ E''_n := \{ x \in \Z_N: T^n \E(f_{U^\perp}|\B'')(x) \geq 3\delta/4 \hbox{ and }
T^n \E(|f_{U^\perp} - f_{UAP}| | \B'')(x) \leq \frac{\delta}{32k} \}.$$
Observe that $E''_n = T^n E''_0$.  Combining this with \eqref{ep-epp}, \eqref{e-ep}, \eqref{mu-squat} we obtain
\begin{equation}\label{mu-est}
 \eqref{mu-form} \geq
 c \E(\prod_{m=1}^{k_*} T^{\mu \lambda m} {\bf 1}_{E''_0}(x) | x \in \Z_N; 1 \leq \lambda \leq N_0/k_*)
- O(N_0^{-30}).
\end{equation}
The function ${\bf 1}_{E''_0}$ is measurable in $\B''$, which is a compact $\sigma$-algebra of order $d-1$.  
At this point it it is tempting to apply the induction hypothesis (Theorem \ref{recurrence} for $d-1$) to
${\bf 1}_{E''_0}$ (using Proposition \ref{compact-sigma}) to obtain lower bounds
for the right-hand side of \eqref{mu-est}.  Unfortunately the problem is that the complexity of $\B''$ depends on
$X'$, whereas the range $N_0/k_*$ of the variable $\lambda$ is only allowed to depend on $X$, and so we cannot
ensure that this expectation is even positive.  To resolve this we must descend from the set $E''_0 \in \B''$
to the slightly modified set $E''_0 \cap \tilde E$, where
$$ \tilde E := \{ \E(f_{U^\perp}|\B) \geq 7\delta/8 \hbox{ and } \E(|f_{U^\perp} - f_{UAP}| | \B) \leq \frac{\delta}{64k}
\}.$$

\begin{lemma}\label{e-bound}  Either we have
$$ \P( \tilde E \backslash E'' ) \leq O(\tau^2); \quad \P( E''_0 \cap \tilde E ) \geq \delta/32$$
or we are in the energy increment half of the dichotomy.  
\end{lemma}

\begin{proof}  We may assume without loss of generality that
$$ \Energy_f(\B'') - \Energy_f(\B') \leq \tau^2$$
since otherwise we would be in the energy increment half of the dichotomy.  From \eqref{egap} we thus have
$$ \Energy_f(\B'') - \Energy_f(\B) \leq 2\tau^2,$$
which implies from \eqref{pythagoras} and definition of $f$ that
$$ \E(  |\E(f_{U^\perp}|\B'') - \E(f_{U^\perp}|\B')|^2 | \Z_N ) \leq 2\tau^2$$
and
$$ \E(  |\E(|f_{U^\perp}-f_{UAP}||\B'') - \E(|f_{U^\perp}-f_{UAP}| |\B')|^2 | \Z_N ) \leq 2\tau^2$$
In particular by Chebyshev's inequality we have
$$ \P( |\E(f_{U^\perp}|\B'') - \E(f_{U^\perp}|\B')| \geq \frac{\delta}{8} ) \leq O(\tau^2)$$
and
$$ \P(  |\E(|f_{U^\perp}-f_{UAP}||\B'') - \E(|f_{U^\perp}-f_{UAP}| |\B')| \geq \frac{\delta}{64k} | \Z_N ) 
\leq O(\tau^2)$$
and the first claim follows from the definitions of $E''_0$ and $\tilde E$.

Now we prove the second claim.  From \eqref{central-1} we have
$$ \E( \E( f_{U^\perp}|\B ) | \Z_N ) = \E(f_{U^\perp} | \Z_N) \geq \delta$$
and hence (by the boundedness of $\E( f_{U^\perp}|\B'' )$)
$$ \P( \E( f_{U^\perp}|\B ) \geq 7\delta/8 | \Z_N ) \geq \delta/8$$
while from \eqref{central-0} and Cauchy-Schwarz we have
$$ \E( \E(|f_{U^\perp} - f_{UAP}| | \B'') | \Z_N ) = \E( |f_{U^\perp} - f_{UAP}| | \Z_N )
\leq \| f_{U^\perp} - f_{UAP} \|_{L^2} \leq \frac{\delta^2}{1024 k}$$
and hence by Chebyshev's inequality
$$ \P( \E(|f_{U^\perp} - f_{UAP}| | \B'') > \frac{\delta}{64 k} | \Z_N ) \leq \delta/16.$$
By definition of $E''_0$, we thus have $\P(E''_0|\Z_N) \geq \delta/16$, and the second claim of the lemma thus follows
from the first if $\tau$ is chosen sufficiently small.
\end{proof}

We may of course assume that we are not in the energy increment half of the dichotomy, in which case
Lemma \ref{e-bound} implies that
$$  \| {\bf 1}_{E''_0 \cap \tilde E} - {\bf 1}_{\tilde E} \|_{L^2} = O(\tau).$$
Now observe that $\tilde E \in \B$ and $\B$ is compact of order $d-1$ and complexity at most $X$.  Thus by Proposition \ref{compact-sigma} we can find a bounded nonnegative function $\tilde f_{UAP} \in UAP^{d-1}$ with
$$ \| {\bf 1}_{\tilde E} - \tilde f_{UAP} \|_{L^2} = O(\tau)$$
and
$$ \| \tilde f_{UAP} \|_{UAP^{d-1}} \leq O_{X,\tau}(1).$$
In view of Lemma \ref{e-bound}, we can thus apply the induction hypothesis of Theorem \ref{e-bound}
with $f_{U^\perp}$ replaced by ${\bf 1}_{E''_0 \cap \tilde E}$, provided that $\tau$ is chosen sufficiently small.  
We conclude that
$$ \E(\prod_{m=1}^{k_*} T^{\mu \lambda m} {\bf 1}_{E_0}(x) | x \in \Z_N; 1 \leq \lambda \leq N_0/k_*)
\geq c(\tau, X)$$
and so if we choose $N_0$ sufficiently large depending on
$\tau, X$, we see
from \eqref{mu-est} that
$$ \eqref{mu-form} \geq c(\tau,X).$$
The claim \eqref{cdm-extend} now follows from \eqref{avast}, and the proof of the Proposition is complete.
\end{proof}

The proof of Szemer\'edi's theorem is now complete.

\section{Appendix: Proof of van der Waerden's theorem}\label{appendix-vwt}

In this section we present the standard ``colour focusing'' proof of van der Waerden's theorem (Theorem \ref{vwt}).  
Our proof presents no new ideas; we give it here only for the sake of 
self-containedness, and also to emphasize that this theorem is comparatively simple compared to
Szemer\'edi's theorem, and thus any argument which manages to reduce the latter to the former is a non-trivial
argument.

The proof of Theorem \ref{vwt} rests on the concept of a \emph{polychromatic fan}, which we now define.
We use $[a,r,k]$ to denote the arithmetic progression $a, a+r, \ldots, a+(k-1)r$.

\begin{definition}  Let $\c: \{1,\ldots,N\} \to \{1,\ldots,m\}$ be a colouring, let $k \geq 1$, $d \geq 0$, and $a \in \{1,\ldots,N\}$.  We define a \emph{fan of radius $k$, degree $d$, and base point $a$} to be a $d$-tuple $([a,r_1,k],\ldots,[a,r_d,k])$ of progressions in $\{1,\ldots,N\}$ of length $k$ and base point $a$, and
refer to the progressions $[a+r_i, r_i, k-1]$, $1 \leq i \leq d$ as the \emph{spokes} of the fan.  We say that a fan
is \emph{polychromatic} if its base point and its $d$ spokes are all monochromatic with distinct colours.   In other words, there exist distinct colours $c_0, c_1, \ldots, c_d \in X$ such that
$\c(a) = c_0$, and $\c(a + j r_i) = c_i$ for all $1 \leq i \leq d$ and $1 \leq j \leq k$.
\end{definition}

\begin{proof}[Proof of Van der Waerden's theorem.]  
We induct on $k$.  The base case $k=1$ is trivial, so suppose $k \geq 2$
and the claim has already been proven for $k-1$; thus for every $m$ there exists a positive integer 
$N_{\vdW}(k-1,m)$ such that any $m$-colouring of $\{1,\ldots,N_\vdW(k-1,m)\}$ contains a monochromatic progression of length $k-1$.

We now claim inductively that for all $d \geq 0$ there exists a positive integer $N_\FAN(k-1,m,d)$ such that
any $m$-colouring of $\{1,\ldots,N_\FAN(k-1,m,d)\}$ contains either a monochromatic progression of length $k$, or
a polychromatic fan of radius $k$ and degree $d$.  The base case $d=0$ is trivial; as soon as
we prove the claim for $d=m$ we are done, as it is impossible in an $m$-colouring for a polychromatic fan
to have degree larger than or equal to $m$.

Assume now that $d > 1$ and the claim has already been proven for $d-1$.  We define $N = N_\FAN(k-1,m,d)$ by the formula
$N := 4k N_1 N_2$, where $N_1 := N_\FAN(k-1,m,d-1)$ and $N_2 := N_\vdW(k-1, m^d N_1^d)$, which are guaranteed to be finite
by the inductive hypotheses, and let $\c$ be an $m$-colouring of
$\{1,\ldots,N\}$.  Then for any $b \in \{1,\ldots,N_2\}$, the set $\{ bkN_1 + 1, \ldots, bkN_1 + N_1\}$
is a subset of $\{1,\ldots,N\}$ of cardinality $N_1$.  Applying the inductive hypothesis, we see that
$\{ bkN_1 + 1, \ldots, bkN_1 + N_1\}$ contains either a monochromatic progression of length $k$, or a polychromatic
fan of radius $k$ and degree $d-1$.  If there is at least one $b$ in which the former case applies, we are 
done, so suppose that the latter case applies for every $b$.  This implies that for every
$b \in \{1,\ldots,N_2\}$ there exist $a(b), r_1(b), \ldots, r_{d-1}(b) \in \{1,\ldots,N_1\}$
and distinct colours $c_0(b), \ldots, c_{d-1}(b) \in \{1,\ldots,m\}$ such that $\c(bkN_1+a(b)) = c_0(b)$
and $\c(bkN_1+a(b)+jr_i(b)) = c_i(b)$ for all $1 \leq j \leq k-1$ and $1 \leq i \leq d-1$.  In particular
the map $b \mapsto (a(b), r_1(b), \ldots, r_{d-1}(b),c_0(b), \ldots, c_{d-1}(b))$ is a colouring
of $\{1,\ldots,N_2\}$ by $m^d N_1^d$ colours (which we may enumerate as $\{1,\ldots,m^d N_1^d\}$ in some arbitrary
fashion).  Thus by definition of $N_2$ there exists a monochromatic
arithmetic progression $[b,s,k-1]$ of length $k-1$ in $\{1,\ldots,N_2\}$, with some colour
$(a,r_1,\ldots,r_{d-1},c_0,\ldots,c_{d-1})$.  We may assume without loss of generality
that $s$ is negative since we can simply reverse the progression if $s$ is positive.  

Now we use an algebraic
trick (similar to Cantor's famous diagonalization trick) which will convert a progression of identical fans into a new fan of one higher degree, the base points of the original fans being used to form the additional spoke of the new fan.  
Introduce the base point $b_0 := (b-s)kN_1+a$, which lies in $\{1,\ldots,N\}$ by construction of $N$, and consider the fan
$$ ([b_0, skN_1,k], [b_0, skN_1 + r_1, k], \ldots, [b_0, skN_1 + r_{d-1}, k])$$
of radius $k$, degree $d$, and base point $b_0$.  We observe that all the spokes of this fan are monochromatic.  For
the first spoke this is because 
$$\c( b_0 + j sk N_1 ) = \c( (b + (j-1) s)kN_1 + a ) = c_0(b + (j-1)s) = c_0 \hbox{ for all } 1 \leq j \leq k-1$$
and for the remaining spokes this is because
$$\c( b_0 + j (sk N_1 + r_t) ) = \c( (b + (j-1) s)kN_1 + a + j r_t ) = c_t(b + (j-1)s) = c_t 
\hbox{ for all } 1 \leq j \leq k-1, 1 \leq t \leq d-1.$$
If the base point $b_0$ has the same colour as one of the spokes, then we have found a monochromatic progression of length
$k$; if the base point $b_0$ has distinct colour to all of the spokes, we have found a polychromatic fan of radius $k$
and degree $d$.  In either case we have verified the inductive claim, and the proof of Proposition \ref{vwt} is complete.
\end{proof}

\begin{remark}  The bounds on $N_\vdW(k,m)$ obtained by this method are clearly of Ackermann type, and are extremely far from best possible.  The first primitive recursive bound on $N_\vdW(k,m)$ is due to Shelah \cite{shelah}.  Currently the best known bound for $N_\vdW(k,m)$ is
$$ N_\vdW(k,m) \leq 2^{2^{m^{c_k}}}, \hbox{ where } c_k := 2^{2^{k+9}},$$
due to Gowers \cite{gowers}.  In contrast to arguments such as the one presented here, in which one deduces
Szemer\'edi-type theorems from van der Waerden type theorems, the bound obtained by Gowers is in fact derived by
the converse procedure, in which one first proves Szemer\'edi's theorem (without recourse to the van der Waerden theorem)
and then deduces the above bound on van der Waerden's theorem as a consequence.
\end{remark}

\begin{remark}  One can modify the above colour focusing technique to prove polynomial or Hales-Jewett versions of
van der Waerden's theorem, see for instance \cite{walters}.
\end{remark}

\end{document}